\documentclass[reqno,11pt,a4paper]{amsart}

\usepackage{amsmath,graphicx,amssymb}
\usepackage[all]{xy}
\usepackage{epic}
\usepackage{eepic}
\usepackage[colorlinks=true,linkcolor=black,citecolor=black]{hyperref}
\newtheorem{theorem}{Theorem}[section]         
\newtheorem{proposition}[theorem]{Proposition}
\newtheorem{corollary}[theorem]{Corollary}
\newtheorem{lemma}[theorem]{Lemma}

\newtheorem{conj}{Conjecture}

\newtheorem{remark}{Remark}

\numberwithin{equation}{section}

\newcommand{\PSL}{\mathrm{PSL}}
\newcommand{\SL}{\mathrm{SL}}
\renewcommand{\le}{\leqslant}
\renewcommand{\ge}{\geqslant}
\renewcommand{\a}{\mathfrak a}
\newcommand{\pp}{\mathfrak p}
\newcommand{\RRp}{\RR^{\mathrm{prim}}}

\renewcommand{\SS}{\mathcal S}
\newcommand{\tSS}{\widetilde{\mathcal S}}
\newcommand{\Z}{\mathbb Z}
\newcommand{\N}{\mathbb N}
\newcommand{\R}{\mathbb R}

\newcommand{\Q}{\mathbb Q}
\newcommand{\RRR}{\mathfrak R}

\renewcommand{\H}{\mathbb H}

\newcommand{\DD}{\mathcal D}

\newcommand{\RR}{\mathcal R}

\newcommand{\MMM}{\mathfrak M}

\newcommand{\NN}{\mathcal N}
\newcommand{\tNN}{\widetilde{\mathcal N}}

\newcommand{\FF}{\mathcal F}

\newcommand{\re}{\operatorname{Re}}

\newcommand{\Vol}{\operatorname{Vol}}

\newcommand{\SSS}{\mathfrak S}

\newcommand{\wQ}{\widetilde{Q}}
\newcommand{\tot}{\mathrm{tot}}
\newcommand{\wom}{\widetilde{\omega}}
\newcommand{\wg}{\widetilde{\gamma}}
\newcommand{\sgn}{\operatorname{sgn}}

\newcommand{\comment}[1]{}

\begin{document}

%



   \author{Florin P. Boca}
   \address{Department of Mathematics, University of Illinois Urbana-Champaign, Urbana, IL 61801, USA}
   \address{Member of the Institute of Mathematics ``Simion Stoilow" of the Romanian
Academy}
   \email{fboca@illinois.edu}


   \author{Alexandru A. Popa}

   \address{Institute of Mathematics ``Simion Stoilow" of the Romanian
Academy  P.O. Box 1-764, RO-014700 Bucharest, Romania}
\email{Alexandru.Popa@imar.ro}
   
   \author{Alexandru Zaharescu}
   \address{Department of Mathematics, University of Illinois Urbana-Champaign, Urbana, IL 61801, USA}
   \address{Institute of Mathematics ``Simion Stoilow" of the Romanian
Academy  P.O. Box 1-764, RO-014700 Bucharest, Romania}
   \email{zaharesc@illinois.edu}


   \title{Pair correlation of hyperbolic lattice angles}


   \begin{abstract}
   Let $\omega$ be a point in the upper half plane, and let $\Gamma$ be
a discrete, finite covolume subgroup of $\mathrm{PSL}_2(\R)$. We conjecture an explicit
formula for the pair correlation of the angles between geodesic rays of the lattice
$\Gamma\omega$,
intersected with increasingly large balls centered at $\omega$. We prove this conjecture for
$\Gamma=\mathrm{PSL}_2(\Z)$ and $\omega$ an elliptic point.
   \end{abstract}





   \date{March 18, 2014}


   \maketitle


\section{Introduction}\label{sec1}

The statistics of spacings measure the fine structure of sequences
of real numbers, going beyond the classical Weyl
uniform distribution. Originating in work of physicists
on random matrices \cite{Wig,Dys}, spacing statistics are conveniently
expressed as the convergence of certain measures, called level
correlations, and respectively level spacing measures. In the past decades these
notions have received significant attention in many areas of
mathematical physics, analysis, probability, and number theory.
For most sequences of interest it is usually very challenging
to prove the existence and describe the limiting spacing measures,
such as the gap distribution or pair correlation,
even when existence is experimentally predicted.

One class of interesting sequences studied in recent years arises from the angular distribution of lattice points.
In the Euclidean scenery one such question is: for a given point
$\mathbf{\alpha}\in \R^2$, describe the statistics of the increasing sequence of finite sets
\[
\Big\{ \frac{\mathbf{m}+\mathbf{\alpha}}{\vert \mathbf{m}+\mathbf{\alpha} \vert} : \mathbf{m}\in\Z^2 \setminus
\{ -\mathbf{\alpha} \} , \vert \mathbf{m}+\mathbf{\alpha}\vert <R\Big\} \subseteq S^1 ,\quad
\mbox{\rm with}\ R\rightarrow \infty,
\]
representing the directions of points in the affine lattice $\mathbf{\alpha}+\Z^2$ with observer located at the
origin.
When $\mathbf{\alpha}\in \Z^2$ and only primitive lattice points are considered, the limiting gap distribution
and the pair correlation were studied and computed in \cite{BCZ0} and \cite{BZ}. The repulsion
between consecutive Farey fractions leads to the vanishing of the corresponding densities on the interval
$[ 0,\frac{3}{\pi^2}]$. When $\mathbf{\alpha} \notin \Q^2$ the gap distribution of this sequence was proved by Marklof
and Str\" ombergsson \cite{MS}
to coincide with the gap distribution of the sequence $(\sqrt{n}\hspace{-3pt}\mod{1})$. The latter was computed by Elkies and McMullen \cite{EM}, with
effective estimates obtained only very recently by Browning and Vinogradov \cite{BV}, building on work of Str\" ombergsson \cite{Str}.
A thorough analysis of the mixed moments of consecutive gaps has been recently undertaken
by El-Baz, Marklof and Vinogradov \cite{MV}, who showed that the pair correlation is Poissonian for Diophantine $\mathbf{\alpha}$.

This paper is concerned with the hyperbolic situation, where the lattice $\Z^2$ is replaced with a lattice (discrete subgroup of finite covolume) $\Gamma$
in $\mathrm{PSL}_2(\R)$. We consider the angles between geodesic rays
$(\omega\rightarrow \gamma\omega)$ in the upper half plane $\H$,
connecting a fixed point $\omega\in\H$ with the (finitely many) points $\gamma\omega$ in its $\Gamma$-orbit, lying
in increasingly large hyperbolic balls.
These angles are well-known to be uniformly distributed (see, e.g., \cite{Nic}) and their uniform distribution in angular sectors can be made effective
\cite{Bo,BKS,Ch,GN,HZ,Iw,RR,RT}.

A first step in the study of the pair correlation of directions of hyperbolic lattice points was completed
in \cite{BPPZ}, where we treated the case $\Gamma=\mathrm{PSL}_2(\Z)$ and $\omega=i$, establishing a formula for the pair
correlation density $g_2(\xi)$ that involves two terms.
The first term is a series over the set of matrices $M$ with nonnegative entries of an explicit
function of $\xi$ depending only on the Hilbert-Schmidt norm of $M$, while the second term is a finite sum involving volumes of bodies
defined in terms of the triangle transformation introduced in \cite{BCZ}.
In this paper we extend the approach introduced in \cite{BPPZ}, conjecturing an explicit formula for the pair correlation
density $g_2(\xi)$ for arbitrary $\Gamma$ and $\omega$. We are able to prove this formula for the
modular group and $\omega$ an elliptic point. Remarkably, we find that $g_2(\xi)$ equals the diagonal value at $\omega$ of
an explicit automorphic kernel.

To state the results, we introduce some notation and definitions.
Let $\omega=u+iv\in \H$ and let $\Gamma$ be a discrete, finite covolume subgroup of $\mathrm{PSL}_2(\R)$.
For $\gamma\in\Gamma$, define $\|\gamma\|:=v\sqrt{2\cosh d(\omega, \gamma \omega)}$,
where $d(z_1,z_2)$ denotes the hyperbolic distance between two points $z_1,z_2 \in \H$. Let $B_Q^\tot$ be the
number of matrices $\gamma\in\Gamma$ in the ball $\| \gamma\| \leqslant Q$, so that asymptotically $B_Q^\tot\sim
\frac{3}{v^2} Q^2$. Let $\theta_\gamma \in (-\pi,\pi]$ denote the angle between the vertical geodesic $[\omega,u]$ and
the geodesic ray $[\omega,\gamma\omega]$. We are interested in the pair correlation density
\[
\begin{split}
g_2(\xi)=& \frac{dR_2(\xi)}{d \xi},  \text{ where } R_2(\xi)=\lim_{Q\rightarrow \infty}
\frac{\RR_Q(\xi)}{B_Q^\tot}, \text{ and } \\
\RR_Q(\xi) =\# \bigg\{ (\gamma,\gamma^\prime)\in \Gamma^2 : &
 \gamma^\prime \neq \gamma,\| \gamma\| \leqslant Q, \| \gamma'\| \leqslant Q \  ,\ 0\leqslant
\frac{1}{2\pi} \big( \theta_{\gamma^\prime}- \theta_\gamma \big) \leqslant
\frac{\xi}{B_Q^\tot}\bigg\}.
\end{split}
\]

\begin{conj}\label{Conj1} Let $\Gamma$ be a discrete subgroup of $\mathrm{PSL}_2(\R)$ with fundamental domain of
finite area $V_\Gamma$. The pair correlation measure
$R_2(\xi)$ exists on
$[0,\infty)$, and is given by a $C^1$ function expressed as a series of three dimensional
volumes. Its density $g_2$ is given by the formula
\begin{equation}\label{1.1}
g_2 \Big( \frac{\xi}{V_\Gamma} \Big)= \frac{V_\Gamma}{\pi\xi^2} \sum_{M\in \Gamma} f_\xi \big(\ell (M) \big),
\end{equation}
where $\ell(M)=d(\omega,M \omega)$ and $f_\xi (\ell)$ is the continuous function defined for
$\ell\geqslant 0$ and $\xi>0$ by
\begin{equation}\label{1.2}
f_\xi ( \ell)=\begin{cases} \vspace{0.15cm}
\ln \Big( \frac{\cosh \ell+\sinh \ell}{\cosh \ell
+\sqrt{\sinh^2\ell-\xi^2}}
\Big)
	 & \text{ if } \xi \leqslant 2 \sinh\big( \frac{\ell}{2}\big), \\  \vspace{0.15cm}
\ln
\Big(\frac{(\cosh \ell+\sinh \ell)(1+\xi^2)}{(\cosh \ell
+\sqrt{\sinh^2\ell-\xi^2})^2}\Big)
	 & \text{ if } 2 \sinh\big( \frac{\ell}{2}\big) \leqslant \xi \leqslant
\sinh \ell,  \\
\ln  \left(\cosh \ell+\sinh \ell \right) =\ell
	 & \text{ if } \sinh \ell \leqslant \xi.
\end{cases}
\end{equation}
\end{conj}
Since the series above is absolutely convergent, by l'Hospital we also deduce the conjectural formula:
\begin{equation}\label{1.3}
g_2(0)=\frac{V_\Gamma}{\pi}\sum_{M\in \Gamma, \ell(M)>0} \frac{1}{e^{2\ell(M)}-1}.
\end{equation}

For the elliptic points for $\mathrm{PSL}_2(\Z)$ we prove the conjecture, using extra symmetries of
the hyperbolic lattices centered at $i$ and $\rho=e^{\pi i /3}$.

\begin{theorem}\label{Thm1.1} Conjecture \ref{Conj1} and formula \eqref{1.3} hold for $\Gamma=\mathrm{PSL}_2(\Z)$ and
for $\omega$ one of the elliptic points $i$ or
$\rho$ (with $V_\Gamma=\frac{\pi}{3}$).
\end{theorem}

\begin{figure}[ht]
\centering
\includegraphics*[scale=0.3, bb= 32 192 570 600]{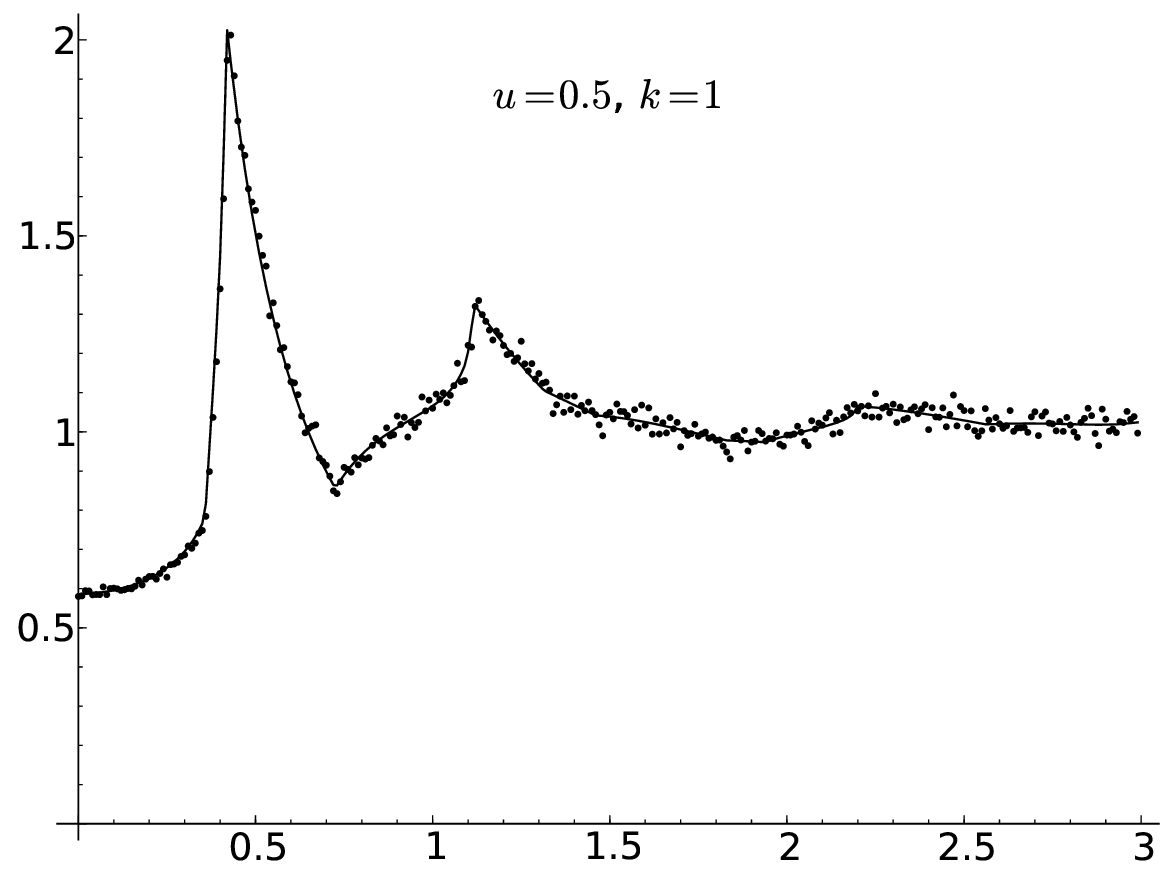} 
\includegraphics*[scale=0.3, bb= 32 192 581 600]{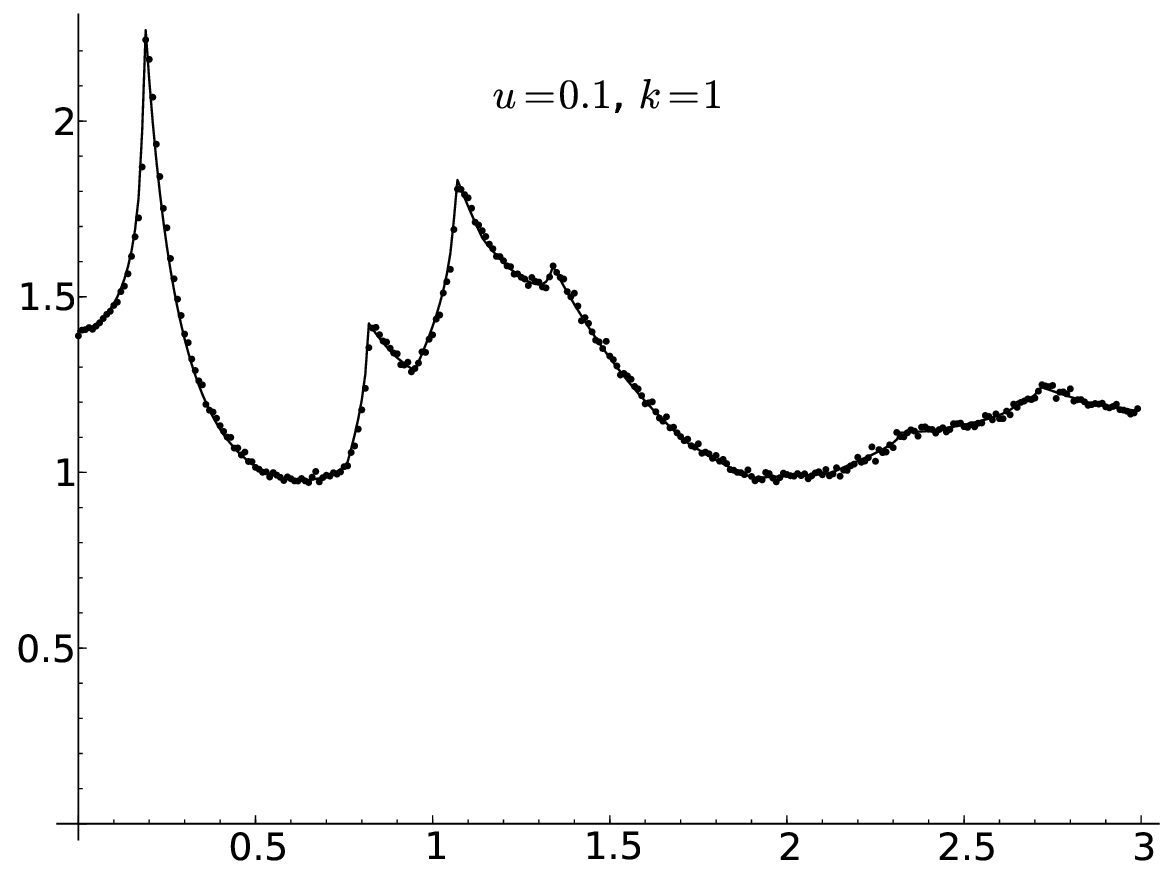} 
\includegraphics*[scale=0.3, bb= 32 192 581 600]{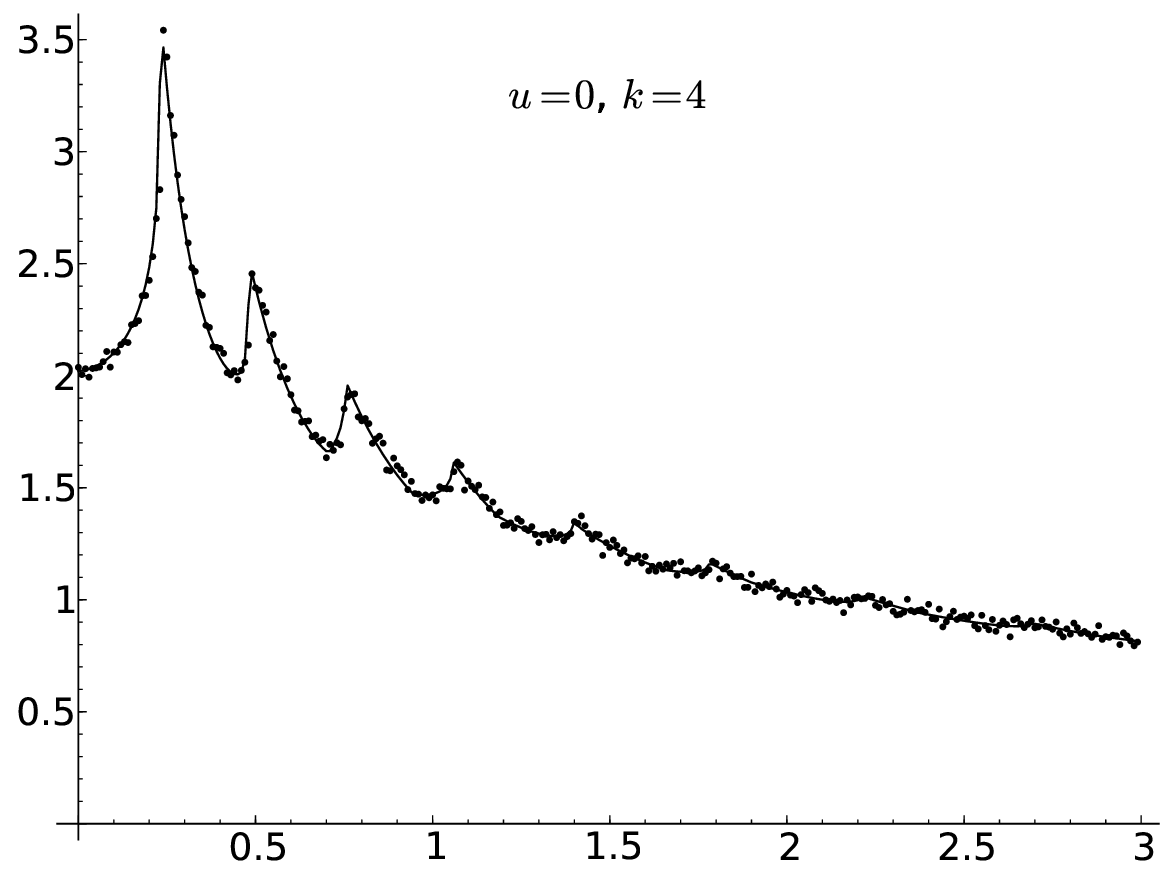}
\includegraphics*[scale=0.3, bb= 32 192 581 600]{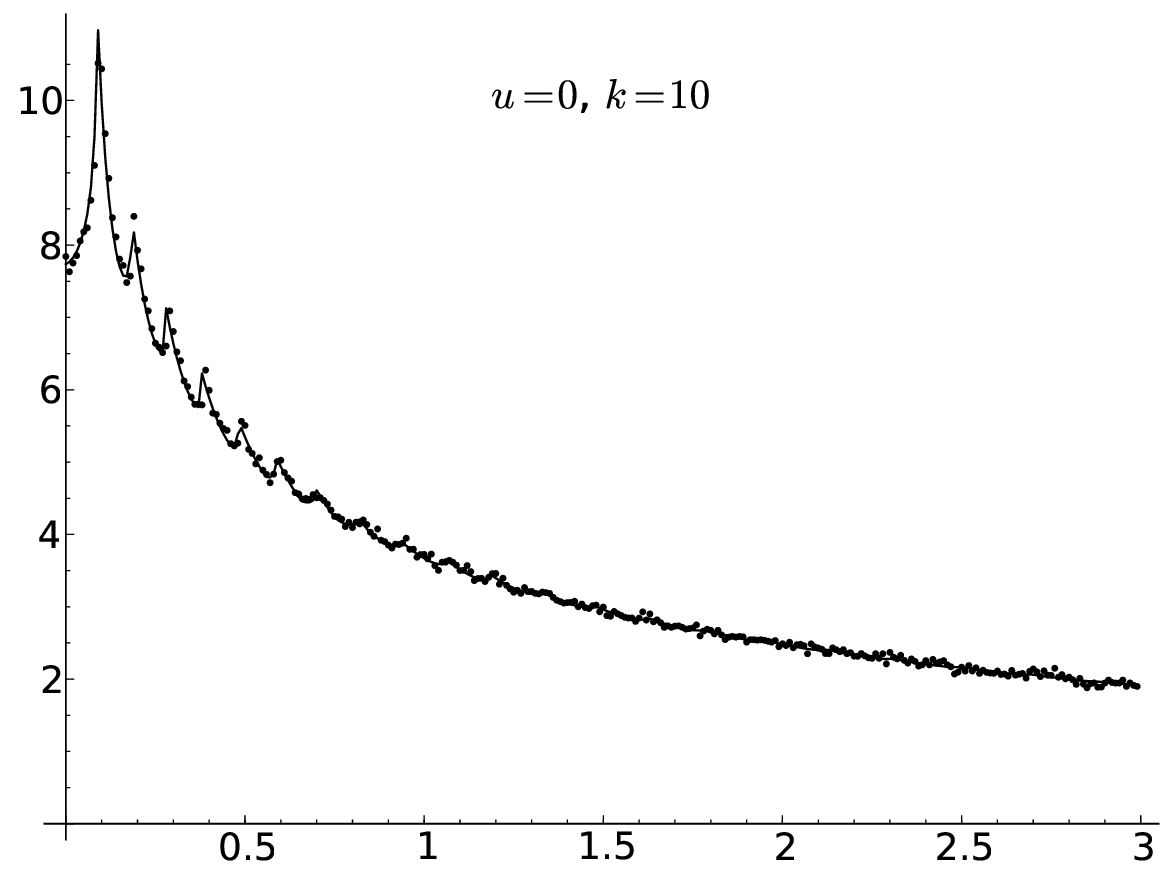}
\caption{The pair correlation functions $g_2^{\mathrm{el}}$ for $\Gamma=\mathrm{PSL}_2(\Z)$ and $\re(\omega)=u$,
$|\omega|=k$ computed using \eqref{1.1} (smooth line) and by counting pairs in the definition
(dots) with $Q=1000$ in the first two plots, $Q=2000$ in the third, and $Q=5000$ in the fourth. }
\label{Figure1}
\end{figure}

The rate of convergence in the result above can be made effective in our proof (see Corollary \ref{C7.5}).
The conjecture has also been verified numerically for a few congruence subgroups $\Gamma_0(N)$ and a few
points $\omega$.
In Fig. \ref{Figure1},  we compare the pair correlation function given
by \eqref{1.1} with the actual pair correlation function computed by counting the pairs in the definition, for
$\Gamma=\mathrm{PSL}_2(\Z)$ and for a few choices of $\omega$. To count the pairs $(\gamma, \gamma')$  in the
definition of $\RR_Q(\xi)$, we first reduce to a half ball $|\gamma \omega|, |\gamma' \omega|\le k$ as explained in Section \ref{sec3}.

When $\Gamma=\mathrm{PSL}_2(\Z)$ and $\omega$ is one of the elliptic points $i$ or respectively $\rho$, each
angle in the definition of $\RR_Q(\xi)$ is counted a number of times equal to the order of
the stabilizer of $\omega$ in $\Gamma$, namely 2 or respectively 3 times. Therefore it is
more natural to consider the pair correlation measure $\RR_Q^{\mathrm{el}}$ defined as $\RR_Q$, with
the condition $\gamma\ne \gamma'$ replaced by $\gamma\omega\ne \gamma'\omega$. Denoting by
$g_2^{\mathrm{el}}$ the corresponding pair correlation functions, we have $g_2^{\mathrm{el}}(\xi)=g_2(e_\omega
\xi)$, where $e_\omega$ is the cardinality of the stabilizer of $\omega$, so that $g_2^{\mathrm{el}}=g_2$ if $\omega$ is
not an elliptic point. For $\omega=i$, the
function $g_2^{\mathrm{el}}$ is identical with the pair correlation function found in \cite{BPPZ},
but the formula here is entirely explicit for all $\xi$.

Formula \eqref{1.1} relates the pair correlation of hyperbolic lattice angles with the length
spectrum of the lattice.  For example, the spikes in the graphs in Fig. \ref{Figure1} occur at values of $\xi$
related in a straightforward way to the length spectrum. Assuming that lattices centered at different
points in a half\footnote{It is shown in Section \ref{sec3} that the pair correlation functions for the lattices
centered at $\omega$ and $-\overline{\omega}$ are equal.} fundamental domain for $\mathrm{PSL}_2 (\Z)$  have different
length spectra, it would follow that the distribution of lattice angles determines the point $\omega$ in a half
fundamental domain.

A common feature of pair correlation density functions, encountered also for the pair correlation of Farey
fractions \cite{BZ}, is that they tend to one at infinity. We expect
the same to be true for the function in Conjecture \ref{Conj1}, i.e. if $\Gamma$ is a Fuchsian group of
the first kind with fundamental domain of finite hyperbolic area $V_\Gamma$, then
\begin{equation}\label{1.4}
\lim_{\xi\rightarrow \infty} \frac{V_\Gamma}{\pi\xi^2} \sum_{M\in \Gamma} f_\xi \big( d(\omega,\gamma\omega) \big) =1,
\end{equation}
where the function $f_\xi(\ell)$ is defined in Conjecture \ref{Conj1}.
The main difficulty in proving this asymptotic formula is that the trace formula does not apply
directly to the automorphic kernel $K_\xi(z,w)=\sum_{\gamma\in\Gamma} k_\xi (u(\gamma z,w))$ with
$k_\xi(u(z,w))=f_\xi (d(z,w))$, because the function $f_\xi$ is not differentiable at
$\xi \in \{ 2\sinh (\frac{\ell}{2}),\sinh \ell\}$.

For $\Gamma=\PSL$ and $\omega=i$, the pair correlation function $g_2^{\mathrm{el}}(\xi)$ is the same as that of
angles made by
reciprocal geodesics on the modular surface, namely the closed geodesics passing through the projection of $i$ on
the modular surface. Reciprocal geodesics were first studied by Fricke and Klein \cite{FK}, and more recently by
Sarnak \cite{Sar}. In \ref{appendix} we similarly describe the arithmetic and geometry of closed geodesics
passing through the projection of $\rho$ on the modular surface. While reciprocal geodesics always consist of two
loops, one tracing the other in the opposite direction, we show that a geodesic on the modular surface passing through the
image of $\rho$ consists of one, two, or four closed loops. The precise
situation depends on the arithmetic properties of the discriminant attached to the geodesic.

We now sketch the main steps in the proof of Theorem \ref{Thm1.1}, while also describing the organization of the
paper. In the remainder of the introduction and throughout the paper we take $\Gamma=\mathrm{PSL}_2(\Z)$, and
keep $\omega$ mostly arbitrary. For technical reasons, we assume that $\re(\omega)$ and $|\omega|^2$ are rational. An
important role is played by the set $\SSS$ of matrices with nonnegative entries, distinct from the identity.

\emph{Step 1.} As in \cite{BPPZ}, our approach is based  on computing the pair correlation of the quantities
$\Psi(\gamma)=u+v\tan (\frac{\theta_\gamma}{2})$ by first approximating them with
$\Phi(\gamma)=\operatorname{Re} (\gamma \omega)$. The reason for preferring the function $\Phi(\gamma)$ is explained
by Lemma \ref{L2.1}, where we show that $\Phi(\gamma)-\Phi(\gamma M)=\Xi_M(c,d)$ with a function $\Xi_M$ depending
only on the lower row $(c,d)$ of $\gamma$.  In Section \ref{sec3} we reduce to angles in the hyperbolic half balls
for which $|\gamma \omega |<|\omega|$, and we show in Section \ref{sec4} that the sets $\{\Psi(\gamma)\}$ and
$\{\Phi(\gamma)\}$ have the same pair correlation.

\emph{Step 2.} To compute the pair correlation of $\{\Phi(\gamma)\}$, we estimate the number $\RR_Q^\Phi(\xi)$ of
pairs $(\gamma, \gamma')\in \Gamma^2$ with $|\gamma \omega|,|\gamma' \omega| <k$, $\|\gamma\|, \|\gamma'\|\leqslant
Q$, and $ 0\leqslant \Phi(\gamma)-\Phi(\gamma')\leqslant\frac{\xi}{Q^2}$ as follows:
\begin{equation}\label{1.5}
\RR_Q^\Phi(\xi)=\frac{1}{2}\sum_{M\in\Gamma\backslash\{I\}} \NN_{M,Q}(\xi) ,
\end{equation}
where  $\NN_{M,Q}(\xi)$ is the cardinality of the set
\begin{equation}\label{1.6}
\SS_{M,Q}(\xi) :=\Big\{ \gamma \in \Gamma: 
\vert \Phi (\gamma M)-\Phi (\gamma)\vert \leqslant \frac{\xi}{Q^2},\
\vert \gamma \omega\vert,\vert \gamma M \omega\vert <k,\  \| \gamma\|,\| \gamma M\| \leqslant Q\Big\} .
\end{equation}
Replacing $b$ in terms of $a,c,d$, for $M\in \SSS$, we show that  $\NN_{M,Q}(\xi)$ is asymptotic as
$Q\rightarrow \infty$ to the cardinality $\tNN_{M,Q} (\xi)$ of the set
$\tSS_{M,Q}(\xi)$ of integer triples $(a,c,d)$ such that
\begin{equation}\label{1.7}
\begin{cases}
\vert a\vert \leqslant k c\leqslant \wQ,\quad \vert d\vert \leqslant \wQ,\quad
ad\equiv 1
\hspace{-6pt} \pmod{c}, \quad
\vert \Xi_M (c,d)\vert \leqslant \frac{\xi}{Q^2},\\  \max\{ c^2 k^2 +d^2 +2cdu,c^2 X_M
+d^2 Y_M +2cd Z_M\} \leqslant \frac{Q^2 c^2}{a^2+k^2 c^2-2acu} ,
\end{cases}
\end{equation}
with $\wQ=\frac{Q\sqrt{k}}{v\sqrt{k-|u|}}$ and $X_M, Y_M, Z_M$ as defined in \eqref{2.1}. This approximation holds for
fixed $M\in \SSS$ with an explicit error term (Lemma \ref{L7.2} (iii)), but in order to control the error
when summing the series \eqref{1.5}, we need to replace $\tNN_{M,Q} (\xi)$ by the cardinality $\tNN_{M,Q}^+ (\xi)$
of the subset $\tSS_{M,Q}^+(\xi)\subseteq \tSS_{M,Q}(\xi)$ consisting of triples $(a,c,d)$ as above with $d>0$.

\emph{Step 3.} Using
estimates for number of points in hyperbolic regions based on bounds on Kloosterman sums, we show in Lemmas
\ref{L7.3} and \ref{L7.4} that for $M\in \SSS$, $\tNN_{M,Q} (\xi)\sim \frac{Q^2}{\zeta(2)} \operatorname{Vol}
(S_{M,\xi})$, and a similar estimate holds also when summing $\tNN_{M,Q}^+ (\xi)$ over $M\in \SSS$. The region
$S_{M,\xi}$ consists of triples $(x,y,z)\in \big[
0,\frac{1}{v\sqrt{k(k-|u|)}}] \times [-\frac{\sqrt{k}}{v\sqrt{k-|u|}},
\frac{\sqrt{k}}{v\sqrt{k-|u|}}] \times [-k,k]$ for which
\begin{equation}\label{1.8}
\begin{cases}
\vert \Xi_M (x,y)\vert \leqslant \xi\quad \mbox{\rm and} \\
\max\{ x^2 k^2 +y^2+2xyu, x^2 X_M + y^2 Y_M +2xy Z_M \} \leqslant \frac{1}{k^2+z^2-2uz} .
\end{cases}
\end{equation}

\emph{Step 4.} Using extra symmetries of the hyperbolic lattice in the case $\omega=i$ (Section \ref{sec8}) and
$\omega=\rho$ (Section \ref{sec10}), we show that the summation range in \eqref{1.5} can be reduced to a subset of
$\SSS$. Moreover, using repulsion arguments involving the Farey tessellation, we can define finite subsets
$\widetilde{\FF}(\xi)\in \Gamma$ such that for $M\in \SSS \backslash \widetilde{\FF}(\xi)$ the quantities
$\tNN_{M,Q}(\xi)$ can be expressed only in terms of $\tNN_{M^\dagger,Q}^+(\xi)$ for appropriate $M^\dagger\in
\SSS$. Therefore we place ourselves in the situations analysed in Steps 2 and 3, obtaining
\begin{equation}\label{1.9}
\RR_Q^\Phi (\xi) \sim \frac{Q^2}{2\zeta(2)} \sum_{M\in \Gamma} \operatorname{Vol} (S_{M,\xi}).
\end{equation}

\emph{Step 5.}  The resulting volumes are expressed in closed form as integrals in Section \ref{sec9} for $\omega$
arbitrary, and passing to the pair correlation function $R_2$ of the angles $\{\theta_\gamma\}$ we obtain the
formula in Conjecture \ref{Conj1}, and finish the proof of Theorem \ref{Thm1.1}.

The difficulty in proving Conjecture \ref{Conj1} for $\Gamma=\mathrm{PSL}_2 (\Z)$ and general $\omega$ resides in the estimates of
Steps 2 and 3, where we  use positivity of some of the entries of the matrices involved to obtain good control on
the error in lattice point counting.

After this paper was completed (arxiv.org/abs/1302.5067), we received a preprint by Kelmer and
Kontorovich \cite{KK} in which they prove Conjecture \ref{Conj1} and the asymptotic formula \eqref{1.4} for a general lattice $\Gamma$ and $\omega\in\H$. Their starting point is a decomposition similar to \eqref{1.5},
for the counting function associated with the angles $\theta_\gamma$ instead of the function $\Phi(\gamma)$, in which they estimate the infinite series by using spectral methods.

\section{Preliminary computations}\label{sec2}

To each $\omega=u+iv \in \H$ with $|\omega|=k$ and $\gamma =\left(\begin{smallmatrix} a & b \\ c & d \end{smallmatrix}\right)\in
\SL_2 (\R)$, we associate
\begin{equation}\label{2.1}
\begin{split}
& X=|a\omega+b|^2, \quad Y=| c\omega+d|^2, \\
& Z=Y\operatorname{Re} (\gamma \omega)= ac|\omega|^2+bd+u(ad+bc),\\
& T =\|\gamma\|^2:= X+|\omega|^2 Y-2u Z .
\end{split}
\end{equation}
Setting $\Delta=2v^2$ and $\epsilon_T=\frac{T-\sqrt{T^2-\Delta^2}}{\Delta}\leqslant \frac{\Delta}{T}$, we have
\begin{equation}\label{2.2}
XY-Z^2=v^2, \ \ \cosh d(\omega, \gamma \omega)=\frac{T}{\Delta},
 \end{equation}
The first equality in \eqref{2.2} leads to
$T=X+k^2 Y-2uZ \geqslant X+k^2 Y -2\vert u\vert \sqrt{XY} =(\sqrt{X}-\vert u\vert \sqrt{Y})^2 +v^2 Y$, and so
\begin{equation}\label{2.3}
v^2 Y \leqslant \| \gamma\|^2 \quad \mbox{\rm and}\quad
X,Y,\vert Z\vert \ll_\omega \| \gamma\|^2 .
\end{equation}

A direct calculation provides
\begin{equation}\label{2.4}
\sin \theta_\gamma=2v\frac{ Z-uY}{\sqrt{T^2-\Delta^2}}, \quad
\cos \theta_\gamma=\frac{\Delta Y-T}{\sqrt{T^2-\Delta^2}}, \quad
\tan( \theta_\gamma/2) =\frac{1}{v}\frac{Z-u Y}{Y-\epsilon_T}.
\end{equation}
The point $\gamma \omega$
is completely determined by the ``coordinates'' $(X,Y,Z)=(X_\gamma,Y_\gamma,Z_\gamma)$, as its hyperbolic polar
coordinates are so determined.

The $x$-intercept $\Psi (\gamma)$ of the oriented geodesic $\omega \rightarrow \gamma\omega$ is given by
\[
\Psi(\gamma)=u+v\tan ( \theta_\gamma / 2 ) =\frac{Z_\gamma -u\epsilon_T}{Y_\gamma-\epsilon_T}.
\]
Since $\epsilon_T=e^{-d(\omega,\gamma \omega)} \rightarrow 0$, a better behaved quantity approximating $\Psi(\gamma)$ well is
\begin{equation*}
\Phi(\gamma):=\frac{Z_\gamma}{Y_\gamma}=\operatorname{Re} (\gamma \omega).
\end{equation*}

For $M=\left( \begin{smallmatrix} A & B \\ C & D \end{smallmatrix} \right)\in
\Gamma$, let $\ell (M)=d(\omega, M\omega)$, and let the angle $\theta_M$ be defined as for
$\gamma$. For $c,d \in \R$, let  $c\omega+d =r e^{i\theta}$,  and define
\begin{equation}\label{2.5}
\Xi_M(c,d) :=-\frac{v}{r^2}
\frac{\sin(\theta_M-2\theta)}{\coth \ell (M)+\cos(\theta_M-2\theta)} .
\end{equation}

\begin{lemma}\label{L2.1}
For $M\in \SL_2(\R)$ and $\gamma=\left( \begin{smallmatrix} a & b \\ c & d
\end{smallmatrix} \right) \in \SL_2(\R)$, we have
\begin{equation*}
\Phi (\gamma)-\Phi (\gamma M) =\Xi_M (c,d).
\end{equation*}
\end{lemma}
\begin{proof}
We compute $\operatorname{Re} (\gamma \omega -\gamma M\omega)$ using the $KAK$ decomposition for $M$ and the
$NAK$ decomposition for $\gamma$, both centered at $\omega$. Denote
$$a(y)=\left( \begin{matrix}
y^{1/2} & 0 \\ 0 & y^{-1/2}
\end{matrix} \right)\in A, \quad n(x)=\left( \begin{matrix} 1 & x \\ 0 & 1
\end{matrix} \right)\in N, \quad k(\theta)=\left( \begin{matrix} \cos \theta &
-\sin \theta \\ \sin \theta & \cos \theta
\end{matrix} \right)\in K
$$
and let $\alpha=n(u)a(v)$, so that $\alpha i=\omega$ and $\alpha$ maps the vertical
through $i$ onto the vertical through $\omega$. Let $z:=\alpha^{-1}\gamma\omega$, so
\begin{equation}\label{2.6}
 \gamma \omega=\alpha z=u+v z , \quad z=x+iy.
\end{equation}
We have
\[\alpha^{-1}\gamma \alpha=  n(x)a(y) k(\theta), \ \ \text{with } c\omega+d =r
e^{i\theta}, r=y^{-1/2}.
\]
Let $\theta_M'$ be the angle between $i\rightarrow i\infty$ and $i\rightarrow \alpha^{-1}
M \alpha i$. Since $\alpha$ maps $i\rightarrow \alpha^{-1} M \alpha i$ onto $\omega
\rightarrow M\omega$, we have
\[
\theta_M'=\pi-\theta_M, \quad d(i,  \alpha^{-1} M \alpha i) = \ell (M).
\]
Consequently
\[
\alpha^{-1}M\alpha=k\big(\theta_M'/2\big) a\big(e^{\ell (M)}\big)k\big(\theta''\big)
\]
for some $\theta''$, and taking $\nu=\theta+\frac{1}{2}\theta_M'$ we obtain
\[
\alpha^{-1} \gamma M \omega = n(x) a(y)k(\theta+\theta_M'/2) a(e^{\ell (M)}) i=x+y
\frac{ie^{\ell (M)}\cos \nu-\sin \nu  }{i e^{\ell (M)} \sin \nu+\cos \nu} .
\]
The equation above and \eqref{2.6} imply
\[
\operatorname{Re} (\gamma M \omega)=\operatorname{Re} (\gamma\omega) +\frac{v}{r^2}\operatorname{Re} \bigg(\frac{ie^{\ell
(M)}\cos \nu-\sin \nu }{i e^{\ell (M)} \sin
\nu+\cos \nu}\bigg).
\]
The real part in the last expression equals the second fraction in \eqref{2.5}, and the claim follows.
\end{proof}

A direction calculation yields
\begin{equation*}
\Xi_M (c,d)= \frac{cd(k^2 Y_M-X_M)+c^2
(k^2 Z_M -uX_M)+d^2(uY_M -Z_M)}{(c^2k^2+d^2+2cdu)(c^2X_M +d^2Y_M +2cdZ_M)}.
\end{equation*}

\section{Reduction to angles with $|\gamma \omega|<|\omega|$}\label{sec3}

We next show that the pair correlation function is determined only by the angles
$\theta_\gamma$ with $|\gamma \omega|<|\omega|=k$, justifying the assumption made in the introduction. When $k\ne 1$,
our proof is conditional upon the formula in Conjecture \ref{Conj1}.

Denote by $\RRR_Q^{\tot}$ the set of $\gamma\in \Gamma$ with $\|\gamma\|<Q$, of
cardinality $B_Q^\tot$, and let $\RRR_Q$, respectively $\RRR_Q^{>}$ denote the subsets for which $|\gamma\omega|<k$,
respectively $|\gamma \omega|>k$, of cardinalities $B_Q$, respectively $B_Q^{>}.$ Lemma \ref{L7.1} will yield
$B_Q\sim B_Q^> \sim \frac{3}{\Delta}Q^2$. Let
\begin{equation*}
\begin{split}
& \RR_Q^{\tot}(\xi) =\# \Big\{ (\gamma,\gamma^\prime)\in (\RRR_Q^{\tot})^2 :
 \gamma^\prime \neq \gamma,\ 0\leqslant
\frac{1}{2\pi} \big( \theta_{\gamma^\prime}- \theta_\gamma \big) \leqslant
\frac{\xi}{B_Q^\tot}\Big\}, \\
& R_2^{\tot}(\xi)=\lim_{Q\rightarrow \infty}\frac{\RR_Q^{\tot}(\xi)}{B_Q^{\tot}},\qquad
g_2^{\tot}=\frac{dR_2^{\tot}}{d\xi},
\end{split}
\end{equation*}
and define similarly $R_2^>$, $g_2^>$. Since $B_Q\sim B_Q^>\sim \frac{1}{2} B_Q^\tot$, we have
$R_2^\tot=\frac{1}{2}(R_2+R_2^>)$,  $g_2^\tot=\frac{1}{2}(g_2+g_2^>)$.

Let $s= \left( \begin{smallmatrix} 0 & -1 \\ 1 & 0 \end{smallmatrix} \right)$. Since
we will compare the hyperbolic lattices centered at the points $\omega$ and $s\omega$, in the
following two paragraphs only we attach subscripts to all notation to denote this
dependence, e.g. $\|\gamma\|_\omega$, $(B_Q)_\omega$ etc.

The map $\gamma\mapsto s \gamma s$ induces the mapping $\gamma \omega \mapsto s\gamma\omega$,
taking the part of the lattice $\Gamma \omega$ with $|\gamma \omega|>|\omega|$,  bijectively and conformally onto the part
of the lattice $\Gamma s\omega$ with $|\gamma s\omega|<|s\omega|$.
Note that
\begin{equation*}
\begin{split}
& k^2 (X_{s\gamma s})_{s\omega}=(X_\gamma)_\omega ,\quad
 k^2 (Y_{s\gamma s})_{s\omega}=(Y_\gamma)_\omega ,\quad
  k^2 (Z_{s\gamma s})_{s\omega}=-(Z_\gamma)_\omega ,\\
  & k^4 \| s\gamma s\|_{s\omega} =\| \gamma\|_\omega ,\quad
  (\theta_{s\gamma s})_{s\omega} =-(\theta_\gamma)_\omega ,
\end{split}
\end{equation*}
yielding
\[(B_Q^>)_\omega=(B_{Q/k^2})_{s\omega},\  \
(\RR_Q^>)_{\omega}=(\RR_{Q/k^2})_{s\omega},\  \   R_2^>(\xi)_\omega = R_2(\xi)_{s\omega}, \  \
g_2^>(\xi)_\omega = g_2(\xi)_{s\omega}.
\]
We conclude that $g_2(\xi)^{\tot}_\omega=\frac{1}{2}( g_2(\xi)_\omega+g_2(\xi)_{s\omega} ).$

Assuming now that $g_2(\xi)_\omega$ is given by the series in Conjecture \ref{Conj1}, we observe that the application
$M\rightarrow s M s$ rearranges the terms of the series for  $g_2(\xi)_\omega$ into the terms of the series for
$g_2(\xi)_{s\omega}$ because the summands only depend on $d(\omega,M\omega)$.  Therefore $g_2(\xi)_\omega=g_2(\xi)_{s\omega}$,
and hence we have $$g_2^\tot(\xi)=g_2(\xi)=g_2^>(\xi)$$ (dropping the subscripts $\omega$ since the basepoint of the
lattice is fixed).

When $k=1$, one can see directly that $g_2(\xi)=g_2^>(\xi)$, because of an extra symmetry of the hyperbolic
lattice. Keeping $\omega$ arbitrary, let $\widetilde{\gamma} =\eta \gamma \eta$ for
$\eta=\left(\begin{smallmatrix} 0 & 1 \\ 1 & 0 \end{smallmatrix}\right)$ and
$\wom=\frac{\omega}{k^2}$. Since
$$
(X_{\wg})_{\wom}=\frac{(Y_\gamma)_{\omega}}{k^2}, \quad
(Y_{\wg})_{\wom}=\frac{(X_\gamma)_{\omega}}{k^2}, \quad
(Z_{\wg})_{\wom}=\frac{(Z_\gamma)_{\omega}}{k^2},
$$
we have $|\gamma \omega|<k \iff |\wg \wom|>k$ and
\begin{equation*}
\big( \Xi_{\widetilde{M}} (y,x)\big)_{\wom}=-\big( \Xi_M(x,y)\big)_{\omega}.
\end{equation*}
The angles $(\theta_\gamma)_\omega$ and
$(\theta_{\widetilde{\gamma}})_{\wom}$ are related as in the following lemma, which shows
directly that $g_2(\xi)=g_2^>(\xi)$ when $|\omega|=1$ (that is $\wom=\omega$).

\begin{lemma}\label{L3.1}
Let $\beta\in (0,\pi)$ be the polar
angle of $\omega=k e^{i\beta}$. We have the relation:
\begin{equation}\label{3.1}
(\theta_\gamma)_\omega+(\theta_{\widetilde{\gamma}})_{\wom}=2 \beta,
\end{equation}
namely the angle between the circle $|z|=k$ and $\omega\rightarrow \gamma \omega$ is the same
as the angle between the circle $|z|=k^{-1}$ and $\wom\rightarrow\wg \wom$.
\end{lemma}

\begin{proof}
When $\omega=ki$ the claim is immediate from \eqref{2.4}. In general, let
$\alpha=\left(\begin{smallmatrix} ak & bk^2 \\ b & ak
\end{smallmatrix}\right)\in \SL_2(\R)$ for
$a=\sqrt{\frac{k+v}{2vk^2}}$, $b=\sqrt{\frac{k-v}{2vk^2}} \sgn u$. Then $\alpha$ fixes the
circle $|z|=k$ and takes $ki$ to $\omega$, while $\widetilde{\alpha}$ fixes the circle
$|z|=k^{-1}$ and takes $k^{-1}i$ to $\wom$, therefore
$$
\angle(\omega\rightarrow
\gamma \omega, \omega\rightarrow k)= \angle( ki \rightarrow
\alpha^{-1} \gamma \alpha ki, k i \rightarrow k).$$
Likewise one has
$$ \angle(\wom\rightarrow
\wg \wom, \wom\rightarrow k^{-1})= \angle( k^{-1}i \rightarrow
\widetilde{\alpha}^{-1} \wg \widetilde{\alpha} k^{-1}i, k^{-1}i \rightarrow k^{-1}),
$$
and the last angles in both equalities are equal by the case $\omega =ki$ already proved.

Alternatively, equality \eqref{3.1} can be checked by direct computation, using the
formula for $\tan(a+b)$.
\end{proof}

From the lemma, combined with $g_2(\xi)_\omega=g_2(\xi)_{s\omega}$ we also deduce that
the pair correlation functions for the hyperbolic lattices centered at $\omega=u+iv$ and
$s\wom=-u+iv$ are equal.  This shows that we can restrict ourselves, whenever convenient, to points $\omega$ in the half
fundamental domain for $\Gamma$ given by
\begin{equation}\label{3.2}
|\omega|\le 1,\quad  \re(\omega)\ge 0, \quad |\omega-1|\ge 1.
\end{equation}

\section{The coincidence of the pair correlations of
$\Phi$ and $\Psi$}\label{sec4}

Since the pair correlation of the lattices centered at $\omega$ and $\gamma_0 \omega$ is the same for $\gamma_0\in
\Gamma=\mathrm{PSL}_2(\Z)$, in this section we assume without loss of generality that $\omega$ lies in a specific fundamental
domain for
the action of $\Gamma$ on the upper half plane. Namely, we assume that $0\leqslant \re (\omega) \leqslant 1$ and
$|\omega-\frac{1}{2}|>\frac{1}{2}$, that is
\begin{equation}\label{4.1}
0\leqslant u\leqslant 1, \quad k^2>u.
\end{equation}
We also need to assume that $u, k^2\in \Q$, which is needed in the proof of Lemma \ref{L4.1}.

Next we show that $\Phi, \Psi$ have the same pair correlation.  As in Section \ref{sec3}, let $\RRR_Q$ be the set
of $\gamma\in \Gamma$ with $|\gamma \omega|<k$ and $\|\gamma\|\leqslant  Q$. Consider
\begin{equation*}
\RR_Q^\Phi(\xi):=\#\Big\{(\gamma,\gamma^\prime)\in (\RRR_Q)^2\ :\
\gamma\ne\gamma^\prime,
\quad 0\leqslant \Phi(\gamma)-\Phi(\gamma')<\frac{\xi}{Q^2} \Big\}
\end{equation*}
and the likewise defined $\RR_Q^\Psi(\xi)$.

Let $\gamma=\left( \begin{smallmatrix} a & b \\ c & d
\end{smallmatrix} \right)\in\Gamma$, $\| \gamma\| \leqslant Q$ and $X,Y,Z,T$
be the quantities defined in the beginning of Section \ref{sec2}. By the results
of Section \ref{sec3} we can restrict to those $\gamma\in \Gamma$ such that $|\gamma
\omega|<k$, that is $X<k^2Y$. In this case we have
$$
\frac{T}{Y}=k^2+\frac{X}{Y}-2u\frac{Z}{Y}<2k^2+2k|u| \ll 1,
$$
and employing $\vert Z-uY\vert \ll T$, a consequence of the first formula \eqref{2.4}, we have
\begin{equation}\label{4.2}
 |\Phi(\gamma)-\Psi(\gamma)|=\frac{\big|u-\frac{Z}{Y} \big|}{\epsilon_T^{-1} Y-1}\ll \frac{1}{Y^2} \ll
\frac{1}{T^2}=\frac{1}{\| \gamma \|^4}.
\end{equation}

Using $XY=Z^2+v^2$ we have $T+2uZ=X+k^2 Y\geqslant 2k\sqrt{XY}>2k|Z|$, hence
$|Z|<\frac{T}{2(k-|u|)}$. It follows
that $\max\{X,k^2Y\}<T+2uZ<\frac{kT}{k-|u|}$.
Since $Y>v^2 \max\{c^2, \frac{d^2}{k^2}\}$, $X>v^2 \max\{ a^2,\frac{b^2}{k^2}\}$, we also
have $\|\gamma\|_\infty \ll \|\gamma\|$, or more precisely:
\begin{equation} \label{4.3}
|a|,|d|<\frac{\|\gamma \|\sqrt{k}}{v\sqrt{k-|u|}},\quad |c|<\frac{\|\gamma
\|}{v\sqrt{k(k-|u|)}}, \quad |b|<\frac{\|\gamma \|k\sqrt{k}}{v\sqrt{k-|u|}}.
\end{equation}

To compare the quantities $\Phi(\gamma)$ and $\Psi(\gamma)$, it will be important to show that they both
lie in a certain Farey interval associate to $\gamma$. For $\Phi(\gamma)$ we can use a
geometric argument to determine this interval. Recall that  $\Phi(\gamma)=\operatorname{Re} ( \gamma \omega)$. Looking
at the images under $\gamma$ of the geodesics $u\rightarrow \omega \rightarrow \infty$ and
$0\rightarrow \omega \rightarrow \frac{k^2}{u}$, if follows that
\begin{equation}\label{4.4}
\Phi(\gamma)\in (\gamma u, \gamma \infty ) \cap \big( \gamma 0,
\gamma (k^2 /u) \big):=J_\gamma^0
\end{equation}
(the endpoints of the intervals are not necessarily ordered increasingly). Since $k^2 Y>X$,
we can assume that $cd(ck^2+du)(cu+d)\ne 0$ at the expense of ignoring a finite number of matrices, which does
not affect the pair correlation. To determine $J_\gamma^0$ explicitly from \eqref{4.4} there are four cases to
consider, and in each one we also define a Farey interval $J_\gamma$ containing $J_\gamma^0$ (using
assumption \eqref{4.1} on $\omega$ and assuming $c>0$):
\begin{enumerate}
 \item[1.] $d>0$. Then $J_\gamma^0=\big(\frac{au+b}{cu+d}, \frac{ak^2+bu}{ck^2+du} \big) \subseteq
\big(\frac{b}{d}, \frac{a}{c}\big)=:J_\gamma$.
 \item[2.] $d<0<cu+d$. Then $ck^2+du>0$, $c+d>0$, $J_\gamma^0=\big(\frac{ak^2+bu}{ck^2+du},
\frac{a}{c}\big) \subseteq \big(\frac{a+b}{c+d},\frac{a}{c}\big)=:J_\gamma$.
\item[3.] $d<0$, $ck^2+du<0$. Then $cu+d<0$, $c+d<0$, $J_\gamma^0 =\big(\frac{b}{d},
\frac{au+b}{cu+d}\big) \subseteq\big(\frac{-b}{-d}, \frac{-(a+b)}{-(c+d)}\big)=:J_\gamma$.
\item[4.] $d<0$, $ck^2+du>0>cu+d$. Then $J_\gamma^0 =\big( \frac{ak^2+bu}{ck^2+du}, \frac{au+b}{cu+d}\big)
\subseteq\big(\frac{a}{c},\frac{-b}{-d}\big)=:J_\gamma$.
\end{enumerate}
With this definition of $J_\gamma$, we have the following   asymptotic
result for $\Psi(\gamma)$.

\begin{lemma}\label{L4.1} Assume $\omega$ satisfies \eqref{4.1} and
let $\gamma=\left( \begin{smallmatrix} a & b \\ c & d
\end{smallmatrix} \right)\in \Gamma$ with $c>0$. Assume also that $u,k^2\in \Q$.
There exists $T_0 =T_0 (\omega)$ such that $\Phi (\gamma),\Psi (\gamma) \in J_\gamma$ whenever $T>T_0$.
\end{lemma}
\begin{proof}
For $\Phi(\gamma)$ the statement was already proved (for all $T$).

Using $cZ -aY =-cu-d$, $dZ-bY=ck^2+du$, we infer
\begin{equation*}
\begin{split}
\Psi(\gamma)-\frac{a}{c}  =\frac{-\frac{cu+d}{c}+\epsilon_T\big(u-\frac{a}{c}\big)}{Y-\epsilon_T} &,\qquad
\Psi(\gamma)-\frac{b}{d} = \frac{\frac{ck^2+du}{d}+\epsilon_T\big(\frac{b}{d}-u\big) }{Y-\epsilon_T}, \\
\Psi(\gamma)-\frac{a+b}{c+d} & =\frac{\frac{c(k^2-u)-d(1-u)}{c+d}+ \epsilon_T\big(\frac{a+b}{c+d}-u\big)
}{Y-\epsilon_T} .
\end{split}
\end{equation*}
We discuss only the second difference, the analysis for the others being similar. We have
$|\frac{ck^2+du}{d}|\gg T^{-1/2}$, since the numerator is bounded from below  as a
result of the rationality assumption on $u$ and $k^2$, and $|d|\le \|\gamma\|_\infty\ll T^{1/2}$.  The
term involving $\epsilon_T$ is $\ll T^{-1}$, thus the sign of $\Psi(\gamma)-\frac{b}{d}$ is the same
as that of $\frac{ck^2+du}{d}$, and similarly for the other two cases, leading to the desired result about $\Psi$.
\end{proof}

\begin{proposition}\label{P4.2}
For each $\beta \in (\frac{1}{2},1)$ one has
\[\RR_Q^\Psi(\xi)<\RR_Q^\Phi(\xi+K_1 Q^{2-4\beta})+ K_2 Q^{1+\beta}\ln Q,
\]
for some constants $K_1, K_2>0$ depending only on $\xi$. The same equality holds with
$\Phi,\Psi$ interchanged.
\end{proposition}

\begin{proof}
Let  $\NN_{Q,\beta}^{>}(\xi)$, respectively $\NN_{Q,\beta}^{<}(\xi)$, be defined as
for $\RR_Q^\Psi(\xi)$, with the additional condition
$\min \{ \| \gamma \|, \| \gamma^\prime \| \}>Q^\beta$, and respectively $\|\gamma\|<Q^\beta$.
We trivially have
\[
\RR_Q^\Psi(\xi)\leqslant\NN_{Q,\beta}^{>}(\xi)+2 \NN_{Q,\beta}^{<}(\xi).
\]
The estimate \eqref{4.2} shows that $\NN_{Q,\beta}^{>}(\xi)\leqslant
\RR_Q^\Phi(\xi+K_1 Q^{2-4\beta})$, the constant $K_1$ being twice the implicit constant in
\eqref{4.2}.

To show $\NN_{Q,\beta}^{<}(\xi)=O_\xi(Q^{1+\beta}\ln Q)$, we follow the same proof
as that of Proposition 3 in \cite{BPPZ}. Because of \eqref{4.3}, at the expense of
counting more pairs we can replace the set $\RRR_Q$ in the definition of
$\NN_{Q,\beta}^{<}(\xi)$ with the set
$$\RRR_Q':=\Big\{\gamma\in\Gamma: |\gamma\omega|<k, \|\gamma\|_\infty\leqslant \wQ :=\frac{Q\sqrt{k}}{v\sqrt{k-\vert
u\vert}} \Big\}.$$
Lemma \ref{L4.1} shows that $\Psi(\gamma), \Phi(\gamma)$ lie between Farey
fractions determined by $\gamma=\left(\begin{smallmatrix} a & b \\ c & d
\end{smallmatrix}\right)$.
More precisely, if $\RRR_Q''$ denotes the subset of $\RRR_Q'$ consisting of matrices
with $cd>0$, and let $I_\gamma=(\frac{b}{d}, \frac{a}{c})$. Then for each $\gamma\in
\RRR_Q'$ we have that
\begin{equation}\label{4.5}
\Psi(\gamma), \Phi(\gamma) \in I_{\gamma'}, \text{ for } \gamma'\in \RRR_Q'',
\end{equation}
with $\gamma'=\gamma$ in Case 1, $\gamma'=\left(\begin{smallmatrix} a & a+b \\ c & c+d
\end{smallmatrix}\right)$ in Case 2, $\gamma'=\left(\begin{smallmatrix} a+b & b \\ c+d & d
\end{smallmatrix}\right)$ in Case 3, and $\gamma'=\left(\begin{smallmatrix} -b & a
\\ -d & c \end{smallmatrix}\right)$  in Case 4 (see the four cases before the statement of the lemma). Note that in
all four cases, $I_{\gamma'}\subseteq [0,k]$ or
$I_{\gamma'}\subseteq[-k,0]$, since $\Phi(\gamma)\in (-k,k)$. Clearly each $\gamma'\in \RRR_Q''$ is associated with one,
two, or three such pairs $\Psi(\gamma),\Phi(\gamma)$ for $\gamma\in\RRR_Q'$.

The proof now follows the same pattern as that of Proposition 3 in \cite{BPPZ},
using \eqref{4.5} above instead of (4.2) there, after further dividing $\RRR_Q''$ into the
subsets of those $\gamma$ with $I_\gamma\subseteq [-k,0]$, and of those $\gamma$ with
$I_\gamma\subseteq [0,k]$. For each subset the analysis of the associated Farey
tessellation formed by the intervals $I_\gamma$ is the same as in  \cite{BPPZ}.
\end{proof}

\section{The Farey tessellation and repulsion}\label{sec5}

In the previous section we associated to each $\gamma \in \Gamma$ an interval $J_\gamma$ between two
consecutive Farey points, such that $\Phi(\gamma)\in J_\gamma$. We also associate to $\gamma$ the geodesic arc on
the upper half plane connecting the endpoints of the Farey interval, which is part of the well known Farey
tessellation. The purpose of
this section is to quantify the statement that there is repulsion between $\Phi(\gamma), \Phi(\gamma')$, if the
intervals associated to $\gamma$, $\gamma'$ are disjoint.\footnote{Recall that two Farey arcs are nonintersecting,
so the corresponding intervals are either disjoint or one contains the other.}

\begin{lemma}\label{L5.1}
Let $\gamma=\left( \begin{smallmatrix} a & b \\ c & d \end{smallmatrix} \right)\in \Gamma$,
$\gamma'=\left( \begin{smallmatrix} e & a \\  f  &  c  \end{smallmatrix} \right) \in \Gamma$ with $ \frac{a}{c}
< \frac{e}{f}$ and $c,d,f >0$. Then there exists $K\in \N$ such that
\begin{equation*}
\gamma \left( \begin{matrix} K & 1 \\ -1 &  0  \end{matrix} \right) = \gamma' .
\end{equation*}
Moreover, if $\max\{Y_{\gamma},Y_{\gamma'}\} \leqslant Q^2$ and we assume $u\ge 0$ and $k\le 1$, then
\begin{equation*}
\Phi (\gamma') - \Phi (\gamma) \geqslant \frac{K k^4}{Q^2} .
\end{equation*}
\end{lemma}

\begin{proof}
Since the matrices $\gamma$ and $\gamma' s$ have the same first column, there exists $K\in \Z$ such that $\gamma
\big( \begin{smallmatrix} 1 & -K \\ 0 &  1  \end{smallmatrix} \big) = \gamma's$, which implies the desired
equality. The fact that $K>0$ follows from $\frac{e}{f}=\frac{aK-b}{cK-d}>\frac{a}{c}$.

A direct calculation provides
\begin{equation*}
\begin{split}
\Phi(\gamma') - \Phi (\gamma) &= \frac{K\alpha (c,d,f)
+cd(1-k^4)+u\big(k^2df+ c^2(1+k^2)-d^2\big)}{(k^2c^2+d^2+2ucd)(k^2f^2+c^2+2ucf) } \\
 &\ge \frac{K \alpha(c,d,f)}{Y_{\gamma}Y_{\gamma'}}\ge \frac{Kk^4}{Q^2} ,
\end{split}
\end{equation*}
where $\alpha(c,d,f)=c^2(k^4+2u^2)+k^2 df+cu(2d+k^2 f)$ denotes the coefficient of $K$ on the first line, and for the first inequality we used
$u\ge 0$, $k\le 1$ and $\frac{c}{d}>K\ge 1 >\frac{1}{\sqrt{1+k^2}}$. If $d\ge f$, then $\alpha(c,d,f)\ge
k^4 Y_{\gamma'}$, and if $f\ge d$, then $\alpha(c,d,f)\ge k^2 Y_{\gamma}$, which, together with
$\max\{Y_{\gamma},Y_{\gamma'}\} \leqslant Q^2$, proves the second inequality.  \end{proof}

For each $\xi >0$ consider the finite set
\begin{equation}\label{5.1}
\begin{split}
\FF(\xi) & :=\bigcup_{\ell\geqslant 1} \bigg\{ M=\gamma_1\cdots \gamma_\ell : \gamma_j =
\left( \begin{matrix} K_j & 1 \\ -1 & 0 \end{matrix}\right),\  K_j\in\N,\   \sum_{j=1}^\ell K_j \leqslant
\frac{4}{k^4}\xi \bigg\}.
\end{split}
\end{equation}

\begin{lemma}\label{L5.2} Assume $u\ge 0, k\le 1$.
Suppose $\gamma=\left( \begin{smallmatrix} a & b \\ c & d \end{smallmatrix} \right),
\gamma^\prime = \left( \begin{smallmatrix} a^\prime & b^\prime \\ c^\prime & d^\prime \end{smallmatrix} \right)
\in \SL_2(\Z)$, $c,d,c^\prime,d^\prime >0$, $\frac{a}{c}\leqslant \frac{b^\prime}{d^\prime}$,
and $\Phi (\gamma^\prime)-\Phi (\gamma) \leqslant \frac{\xi}{Q^2}$ for some $Q\geqslant \max \{
c,d,c^\prime,d^\prime\}$. Then

{\em (i)} $\gamma^\prime=\gamma M$ for some $M\in\FF(\xi)$.

{\em (ii)} Furthermore, if $M=\gamma_1\cdots \gamma_\ell$ is as in \eqref{5.1}, then
$\gamma \gamma_1 \cdots \gamma_j =\left( \begin{smallmatrix} a_j & a_{j-1} \\ q_j & q_{j-1} \end{smallmatrix}\right)$,
with  $q_1,\ldots, q_\ell \in \{ 1,2,\ldots, Q\}$.
\end{lemma}

\begin{proof}
As $Q\geqslant \max\{ c,d^\prime\}$ the fractions $\frac{a}{c}$ and $\frac{b^\prime}{d^\prime}$
belong to the set $\FF_Q$ of ``extended" Farey fractions $\frac{a}{q}$ with
$(a,q)=1$ and $1\leqslant q\leqslant Q$.
Let $\frac{a}{c} =\frac{a_0}{q_0} < \frac{a_1}{q_1} < \cdots < \frac{a_\ell}{q_\ell} =\frac{b^\prime}{d^\prime}$ be
the
elements in $\FF_Q$ between $\frac{a}{c}$ and $\frac{b^\prime}{d^\prime}$. By Lemma \ref{L5.1} there
are positive integers $K_1$ and $K_{\ell+1}$ such that
\begin{equation}\label{5.2}
\left( \begin{matrix} a_1 & a_0 \\ q_1 & q_0 \end{matrix}\right) = \left(
\begin{matrix} a & b \\ c & d \end{matrix}\right) \left( \begin{matrix}
K_1 & 1 \\ -1 & 0 \end{matrix}\right), \quad \ldots\quad
 ,\left( \begin{matrix} a^\prime & b^\prime \\ c^\prime & d^\prime \end{matrix}\right) =
\left( \begin{matrix} a_\ell & a_{\ell-1} \\ q_\ell & q_{\ell-1} \end{matrix}\right)
\left( \begin{matrix} K_{\ell+1} & 1 \\ -1 & 0 \end{matrix} \right) .
\end{equation}
The recursion relations for consecutive Farey fractions
$\frac{q_{j-1}+q_{j+1}}{q_j} =K_j =\frac{a_{j-1}+a_{j+1}}{a_j}$, $j=2,\ldots,\ell-1$, $\ell\geqslant 2$,
\eqref{5.2}, and the consequence $\Phi(\gamma^\prime)-\Phi(\gamma) \geqslant \frac{1}{Q^2}(K_1+\cdots +K_{\ell+1})$
of Lemma \ref{L5.1} yield both (i) and (ii). Note that $Y_{\gamma}, Y_{\gamma'} \le Q^2(k^2+1+2u)\le 4 Q^2$ under
the assumptions on $c,d,c',d', k,u$.
\end{proof}

\section{The case where $\ell(M)$ is large}\label{sec6}

In this section we generalize Lemma 9 of \cite{BPPZ}. Let $\SS^+_{M,Q}(\xi)$, respectively
$\SS^-_{M,Q}(\xi)$,
denote the subsets of the set $\SS_{M,Q}(\xi)$ defined in \eqref{1.6} consisting of
matrices with $c,d>0$, respectively $c>0>d$. Denote by $\NN^\pm_{M,Q}(\xi)$ the cardinality of
$\SS^\pm_{M,Q}(\xi)$.

\begin{lemma}\label{L6.1}
Assume $\omega$ is in the half fundamental domain given by \eqref{3.2}, and $u,k^2\in \Q$.
Suppose $\beta_0 \in ( \frac{1}{2},1)$, $M\in \Gamma$ has nonnegative entries and
$\max\{ X_M,Y_M\} \geqslant Q^{2\beta_0}$.
There exists $Q_0=Q_0(\xi,\omega)$
independent on $M$ such that $\NN^+_{M,Q}(\xi)=0$ for $Q\geqslant Q_0$.
\end{lemma}

\begin{proof}
We show that the region $\Omega_{M,Q}(\xi)$ of $(c,d)\in (0,\infty)^2$ for which
\begin{equation*}
\vert \Xi_M (c,d)\vert\leqslant \frac{\xi}{Q^2} ,\qquad
v^2 \max\{ Y_{\gamma},Y_{\gamma M} \} \leqslant Q^2 ,
\end{equation*}
contains no coprime integer lattice points. Because of \eqref{2.3} this gives $\NN^+_{M,Q}(\xi)=0$.

Suppose there is $(c,d) \in \Omega_{M,Q} (\xi)\cap\Z^2$.
Write $c\omega +d=re^{i\theta}$, and let $X=X_M$, $Y=Y_M$, $Z=Z_M$ and $T=T_M$ be given by \eqref{2.1}.
With $U_M=\coth \ell (M)=1+O( \frac{1}{T^2})$, the inequalities in the definition of $\Omega_{M,Q}(\xi)$ can be described as
\begin{equation}\label{6.1}
\frac{v}{\xi} \frac{\vert\sin (\theta_M -2\theta)\vert}{U_M +\cos (\theta_M -2\theta)} \leqslant
\frac{r^2}{Q^2} 
\leqslant \min\left\{ \frac{1}{v^2},\frac{2}{\sqrt{T^2-\Delta^2} (U_M +\cos (\theta_M
-2\theta))}\right\} .
\end{equation}

Since $\sin \theta>0$, $\cos \theta>0 $ we can take $\theta \in ( 0,\frac{\pi}{2})$. Denoting
$\delta_M=\frac{\theta_M}{2}-\theta $, from  the first and
last fraction in \eqref{6.1} we have
$\vert \sin  2 \delta_M \vert \ll \frac{1}{T}$. Therefore $\delta_M$ is close to 0, or to
$\pm \frac{\pi}{2}$. When $\delta_M$ is close to $0$ we have
$|\tan \delta_M|\ll \vert \delta_M\vert \ll  \vert \sin 2\delta_M \vert \ll\frac{1}{T}$.

When $\delta_M$ is close to $\pm \frac{\pi}{2}$ we similarly have
$|\delta_M\mp \frac{\pi}{2}|\ll \frac{1}{T},$ which we claim is impossible. To prove this, we will use
\begin{equation}\label{6.2}
 \frac{|\tan\delta_M|}{1+\frac{U_M-1}{1+\cos 2\delta_M}}
=\frac{\vert \sin 2\delta_M \vert}{U_M+\cos 2\delta_M} \leqslant \frac{2\xi}{v^3},
\end{equation}
analyzing the two cases that can occur. Note that the equality $Z=uY$ would give $v^2=Y(X-u)$, which cannot hold.

\noindent{\bf Case I:} $Z>uY$. In this case $\theta_M \in (0,\pi)$. Since $u,k^2\in \Q$, \eqref{2.4} yields $ \sin \theta_M\gg \frac{1}{T}$.
Using $\sin\theta_M \sin\theta \geqslant 0$ we infer
\[\textstyle
1+\cos 2\delta_M \geqslant 1+\cos\theta_M \cos\theta \geqslant 1 -\vert \cos
\theta_M \vert = 1-\sqrt{1-\sin^2 \theta_M} \gg \frac{1}{T^2} .
\]
Since $U_M-1\ll \frac{1}{T^2}$, this gives  $0<\frac{U_M-1}{1+\cos 2\delta_M} \ll 1$.
We also have $|\delta_M \pm \frac{\pi}{2}| \ll \frac{1}{T}$, hence $| \tan\delta_M| \gg T$ and the left-hand side in \eqref{6.2}
becomes $\gg T$, producing a contradiction when $Q$ is chosen large enough.

\noindent{\bf Case II:} $Z<uY$. In this case $\theta_M \in (-\frac{\pi}{2},0)$ because $\sin\theta_M <0$ and
$\cos\theta_M>0$ as a result of \eqref{2.4} and of
$$
\Delta Y -T=2uZ -X +(v^2-u^2)Y > uZ +(v^2-u^2)Y -v^2/Y > 0
$$
for $Q$ large enough, where we used $X Y -Z^2=v^2$, and $v>u$ for $\omega$ in the region defined by \eqref{3.2}.
Hence $\delta_M \in ( -\frac{3\pi}{4},0)$, and we will show that $|\delta_M+ \frac{\pi}{2}|\ll
\frac{1}{T}$ leads to contradiction.

As $\vert 2\delta_M \pm\pi\vert \ll \frac{1}{T}$, from $\vert \sin 2\delta_M\vert \ll \frac{1}{T}$ we have
$\vert \tan\delta_M\vert \gg T$, and \eqref{6.2} gives $1\gg  T \frac{1+\cos 2\delta_M}{U_M+\cos 2\delta_M}$. This leads
(for large $T$) to $\frac{1}{T^2} \gg U_M -1 \gg T(1+\cos 2\delta_M)$, and therefore to
$U_M+\cos 2\delta_M \ll \frac{1}{T^2}$. Back to \eqref{6.2}, we infer
$\vert \sin 2\delta_M \vert \ll \frac{1}{T^2}$, and thus
\begin{equation*}
T^2 \ll \vert \tan\delta_M\vert =\frac{\big| \tan \big( \frac{\theta_M}{2}\big) -\frac{cv}{cu+d}\big|}{\big| 1+
\frac{cv}{cu+d} \tan \big( \frac{\theta_M}{2} \big)\big|} <
\frac{1+\frac{v}{u}}{\big| 1+ \frac{cv}{cu+d} \tan \big(\frac{\theta_M}{2}\big) \big|} ,
\end{equation*}
where we used $\theta_M\in ( -\frac{\pi}{2},0) $. This further gives
\begin{equation*}
\bigg| 1+\frac{cv}{cu+d} \tan \big( \theta_M /2 \big) \bigg| =\bigg|
1+\frac{cv}{cu+d} \bigg( \frac{1}{v}\Big( \frac{Z}{Y}-u\Big) +O\Big( \frac{1}{TY}\Big)\bigg)\bigg| \ll_\omega \frac{1}{T^2} ,
\end{equation*}
and so, employing also \eqref{2.3}, $X\ll_\omega Y$ and $Y\gg_\omega Q^{2\beta_0}$, we find
\[
\bigg| 1+\frac{c}{cu+d} \bigg( \frac{Z}{Y}-u\bigg)\bigg| \ll_\omega \frac{cv}{cu+d}
\frac{1}{TY} +\frac{1}{T^2} \leqslant  \frac{v}{uTY}+\frac{1}{T^2} \ll_\omega \frac{1}{TY}\ll \frac{1}{Q^{2\beta_0}Y} .
\]
Finally, multiplying by $\frac{cu+d}{c}$ we infer (here $Z\geqslant ACk^2+BD\geqslant k^2$ and $Y\geqslant k^2$) that there exists a constant
$K=K(\omega,\beta_0)>0$ such that
\[
0< \frac{d}{c} +\frac{Z}{Y} \leqslant \bigg( 1+\frac{d}{c}\bigg) \frac{Kk^2}{Q^{2\beta_0} Y}
\leqslant \frac{d}{c} \frac{K}{Q^{2\beta_0}} +\frac{1}{Y} \frac{Kk^2}{Q^{2\beta_0}}
\leqslant \bigg( \frac{d}{c}+\frac{Z}{Y}\bigg) \frac{K}{Q^{2\beta_0}} ,
\]
which gives a contradiction for $Q$ large.

We have thus shown that $|\tan \delta_M|\ll\frac{1}{T}$. Note that integrality of $c,d>0$ was not used in this argument.

Next we consider the two cases $\tan \theta<v$ and $\tan \theta>v$.

{\bf Case (A).} $\tan \theta < v$. Recall that $\Psi(M)=u+v \tan(\frac{\theta_M}{2})$, and we also have
$u+v \tan
\theta=\frac{ck^2+du}{cu+d}$. Since $\vert \delta_M\vert \ll \frac{1}{T}$, $\tan (\frac{\theta_M}{2})$ is also
bounded, leading to
\begin{equation*}
\begin{split}
\left|\Psi(M) -\frac{ck^2+du}{cu+d}\right| & =v\left| \tan \big ( \theta_M /2 \big)-\tan
\theta \right| \\  &
=v \vert \tan \delta_M \vert \left| 1+\tan\theta \tan \big( \theta_M /2 \big)\right|
\ll Q^{-2\beta_0} .
\end{split}
\end{equation*}
Since $\Psi(M)\ll 1$, we have $Z\ll Y$, and from $X Y -Z^2=v^2$ we conclude $X\ll Y$.
From \eqref{4.2} it follows that $|\Psi(M)-\Phi(M)|\ll Y^{-2}\ll Q^{-4\beta_0}$, so
$ \big| \frac{Z}{Y} -\frac{ck^2+du}{cu+d}\big|   \ll  Q^{-2\beta_0}$. On the other hand,
$$
\frac{A}{C}-\frac{Z}{Y}=\frac{D+uC}{CY} \ll \frac{1}{Y}\ll Q^{-2\beta_0} ,
$$
where we assumed without loss of generality $C\ge D$. If $D\ge C$, then use $\frac{B}{D}$ instead  of $\frac{A}{C}$.
We conclude that
\begin{equation}\label{6.3}
 \left| \frac{A}{C}-\frac{ck^2+du}{cu+d}\right|\ll Q^{-2\beta_0}.
\end{equation}
If nonzero, the left hand side of \eqref{6.3} is $\gg_\omega
\frac{1}{C(cu+d)}$, using the rationality assumption $u,k^2\in \Q$. From $\tan \theta <v$ and $u<1$ it follows that
$c\ll d$, so $Q^{2\beta_0} \ll_\omega C(cu+d)\ll \sqrt{d^2 Y} <\sqrt{Y_{\gamma M}} \ll Q$, which gives a
contradiction.  It remains that $\frac{A}{C}=\frac{ck^2+du}{cu+d}=\frac{Mc+Nd}{Pc+Rd} $ with $M,N,P,R\in \N$
constants, so $d \gg_\omega Pc+Rd\ge C\gg \sqrt{Y}$ (using $C\ge D$). It follows that $Q^2>Y_{\gamma M}>d^2 Y\gg Y^2 \gg
Q^{4\beta_0}$, which again provides a contradiction.

\noindent{\bf Case (B).} $\tan \theta > v$. We have $\frac{u}{v}<\cot \theta<\frac{1}{v} $, and  since
$|\delta_M|\ll \frac{1}{T}$ it follows that $\cot (\frac{\theta_M}{2})$ is bounded as well. Consequently
$$
\left|\frac{1}{\Psi(M)} -\frac{cu+d}{ck^2+du}\right| =v\vert \tan \delta_M \vert \frac{|1+\cot \theta \cot
(\theta_M/2)|}{|(u\cot \theta+v)(u\cot (\theta_M/2)+v)|}\ll Q^{-2\beta_0}.
$$

In this case $\Psi(M)\gg 1$ implies $Z\gg Y$, so $X\gg Z\gg Y$, and from $XY-Z^2=v^2$ we have
$Z\gg Q^{\beta_0}$. Taking into account \eqref{2.4} we arrive at
$$
\left|\frac{1}{\Psi(M)} -\frac{1}{\Phi(M)}\right|=\epsilon_T\frac{\vert
1-u\frac{Y}{Z}\vert}{Z-u\epsilon_T} \ll \frac{1}{T Z}\ll Q^{-3\beta_0}.
$$
Assuming $A>B$, we infer
\[
\frac{Y}{Z}-\frac{C}{A}=\frac{D+uC}{AZ} \ll \frac{\sqrt{Y}}{\sqrt{X}Z}\ll \frac{1}{X}\ll Q^{-2\beta_0} ,
\]
and thus $| \frac{C}{A}-\frac{cu+d}{ck^2+du}| \ll Q^{-2\beta_0}$, leading to a contradiction
for large $Q$ as before. The case $A<B$ is similar, replacing $\frac{C}{A}$ by $\frac{D}{B}$.
\end{proof}

\section{Approximating the number of lattice points in planar regions by volumes}\label{sec7}

In this section we approximate $\NN_{M,Q}(\xi)$ for $M\in{\mathfrak S}$ by volumes of three dimensional regions,
where ${\mathfrak S}\subseteq \Gamma$ is the set of matrices with nonnegative entries, distinct from the identity.
Using the result of the previous section, we also show that the sum of $\NN_{M,Q}^+(\xi)$
over subsets of $M\in{\mathfrak S}$ can be approximated by the corresponding sum of volumes.

To count points in two dimensional regions we use Lemma 7 of \cite{BPPZ}. The prototype for its application in
the present setting is given in the following simpler counting problem. By well known asymptotics for the number of
points in expanding hyperbolic balls we have $B_Q^\tot\sim \frac{6Q^2}{\Delta}$ (for the notation see Section
\ref{sec3}). In the next lemma we show that in half balls we have half this number of points.

\begin{lemma}\label{L7.1}
Let $B_Q=\# \RRR_Q$ be as defined in Section \ref{sec3}, with $k\geqslant 1$. Then
\[B_Q=\frac{3Q^2}{\Delta}+O_\varepsilon (Q^{11/6+\varepsilon}),
\]
and so $B_Q\sim\frac{1}{2} B_Q^\tot.$
\end{lemma}

\begin{proof}
Replacing $b=\frac{ad-1}{c}$, the condition $\vert \gamma \omega\vert <k$ is
equivalent to
\begin{equation}\label{7.1}
\Big(k^2-\frac{a^2}{c^2}\Big) (c^2 k^2+d^2+2cdu)+\frac{2ad-1}{c^2}+\frac{2au}{c} >0.
\end{equation}
A direct calculation shows that $\| \gamma\|\leqslant Q$ is equivalent to
\begin{equation}\label{7.2}
\Big( k^2+\frac{a^2}{c^2} -\frac{2au}{c}\Big) (c^2 k^2 +d^2 +2cdu)
+ \frac{2du}{c}+\frac{1-2ad}{c^2} +\Big( 2u-\frac{2a}{c}\Big) u \leqslant Q^2.
\end{equation}

Fix $\alpha=\frac{13}{18}$ and let $\gamma=\left(\begin{smallmatrix} a & b \\ c
& d \end{smallmatrix}\right)\in \RRR_Q$ with $c>0$. The contribution of the matrices $\gamma$ with
$|c|<Q^{\alpha}$ or $|d|<Q^\alpha$ to the error term is $\ll Q^{1+\alpha}$, so we can assume
$c>Q^{\alpha}$ and $|d|>Q^\alpha$.

We show that the matrices $\gamma$ with $\frac{a^2}{c^2}\geqslant k^2$ contribute negligibly to
$B_Q$. From \eqref{7.1} and $c^2k^2+d^2+2cdu \geqslant c^2 v^2$, it follows that for such $\gamma$ we have
$0\leqslant \frac{a^2}{c^2}-k^2\ll \frac{1}{c^2} ( \frac{Q^2}{c^2}+\frac{Q}{c}) \ll Q^{2-4\alpha}$,
so $\frac{\vert a\vert}{c}\in [k,k+m]$ with $m\ll \sqrt{k^2+K Q^{2-4\alpha}} -k
\ll Q^{1-2\alpha}$ for $K>0$ fixed constant. Since $|a|,c\ll Q$,
from the equidistribution of the Farey
fractions $\mathcal{F}_Q$ in intervals $I$ of length $|I|\gg Q^{-\delta}$ with $\delta=2\alpha -1 \in
(0,1)$ it follows that the number of pairs $(a,c)$ is $\ll Q^{3-2\alpha}$ as long as $\alpha\in ( \frac{1}{2},1)$.
Since the number of values $d$ can
take is $\ll \frac{Q}{c} < Q^{1-\alpha}$, there are $\ll Q^{4-3\alpha}=Q^{11/6}$ such
matrices, so they can be absorbed in the error term.

Therefore we can assume $|a| < k c$, and the condition $|\gamma \omega|<k$ is satisfied
except for a negligible number of matrices. Via \eqref{7.2}, the condition
$\|\gamma\|\leqslant Q$ can be replaced without affecting the asymptotics by
\[
\Big( k^2+\frac{a^2}{c^2}-\frac{2au}{c}\Big)(k^2 c^2+d^2+2cdu)\leqslant Q^2 .
\]
Therefore we can apply Lemma 7 in \cite{BPPZ}, with $q=c$, $L=c^{5/6}$, to the set
$$
\Omega_c=\left\{ (a,d)\in [-kc,kc]\times [-\wQ,\wQ]: (k^2c^2+a^2-2uca)(k^2c^2+d^2+2ucd)\leqslant
Q^2 c^2\right\}
$$
(recall $\wQ=\frac{Q\sqrt{k}}{v\sqrt{(k-|u|)}}$) with area $A (\Omega_c) \ll cQ$ and boundary length
$\ell(\partial \Omega_c) \ll Q$, and conclude that
$$
B_Q=\sum_{c=1}^{\wQ/k}
\frac{\varphi(c)}{c}\frac{A(\Omega_c)}{c}+O_\varepsilon (Q^{11/6+\varepsilon}).
$$
M\"obius summation and a change of variables $a=cz, d=Qy, c=Qx$ then gives
$$
B_Q=\frac{Q^2}{\zeta(2)}\Vol(V_Q)+ O_\varepsilon (Q^{11/6+\varepsilon}),
$$
with $V_Q$ denoting the set of triples $(x,y,z)\in [ 0,\frac{\wQ}{kQ} ]\times [ -\frac{\wQ}{Q},\frac{\wQ}{Q}] \times[-k,k]$ such that
$(k^2x^2+y^2+2uxy)(k^2+z^2-2uz)\leqslant 1$.
The substitution $x\omega+y=re^{i\theta}$, $z=v\tan t+u$, then yields
\[
\Vol(V_Q)=\int_{\beta/2-\pi/2}^{\beta/2}\int_{0}^{\pi} \int_0^{\frac{\cos t}{v}}
\frac{r}{\cos^2 t} dr d\theta dt=\frac{\pi^2}{2\Delta},
\]
which concludes the proof.
\end{proof}

Next we seek to replace inequalities defining the set
$\SS_{M,Q}(\xi)$  in \eqref{1.6}
with simpler ones, involving only the entries $(a,c,d)$ or $(b,c,d)$ of the matrix
$\gamma=\left(\begin{smallmatrix} a & b \\ c
& d \end{smallmatrix}\right)$. Using
\begin{equation*}
\begin{split}
& X_{\gamma M} =a^2 X_M +b^2 Y_M +2ab Z_M, \quad
Y_{\gamma M} =c^2 X_M +d^2 Y_M +2cd Z_M, \\
& Z_{\gamma M}  = acX_M +bdY_M +(ad+bc) Z_M ,
\end{split}
\end{equation*}
and substituting $b=\frac{ad-1}{c}$, we find
\begin{equation}\label{7.3}
\| \gamma M\|^2 = \Big( k^2+\frac{a^2}{c^2} -2u\frac{a}{c} \Big) Y_{\gamma M}
+\frac{Y_M+2(uc-a)(dY_M+cZ_M)}{c^2} ,
\end{equation}
\begin{equation*}
\vert \gamma M \omega\vert <k\ \Longleftrightarrow \
\Big( k^2-\frac{a^2}{c^2}\Big) Y_{\gamma M} +\frac{2a(dY_M+cZ_M)-Y_M}{c^2} >0 .
\end{equation*}

The previous formulas lead us to consider the cardinality
$\tNN_{M,Q} (\xi)$ of the set $\tSS_{M,Q}(\xi)$ of integer triples $(a,c,d)$ satisfying \eqref{1.7}.
Let $\tSS^+_{M,Q}(\xi)$, respectively $\tSS^-_{M,Q} (\xi)$, be the subsets of
$\tSS_{M,Q}(\xi)$ for which  $c,d>0$, respectively $c>0>d$.

\begin{lemma}\label{L7.2} Assume $u,k^2 \in \Q$.

{\em (i)} There is a constant $K=K(\xi)>0$ such that, for every $Q$, the
number of pairs $(\gamma,\gamma^\prime)\in \Gamma^2$ with $\|\gamma\|,\|\gamma^\prime\|\leqslant Q$,
$\min \{\vert \gamma \omega\vert, \vert \gamma^\prime \omega \vert\}< k < \frac{|a|}{c}$ or $\max \{\vert \gamma
\omega\vert, \vert \gamma^\prime \omega \vert\}> k > \frac{|a|}{c}$, and
\begin{equation}\label{7.70}
 \vert \Phi (\gamma^\prime )-\Phi (\gamma)\vert \leqslant \frac{\xi}{Q^2} ,
\end{equation}
is at most $K$.

{\em (ii)} For each $\alpha \in (0,1)$ the following asymptotic estimates hold:
\begin{equation}\label{7.5}
\sum_{M\in {\mathfrak S}} \NN_{M,Q}^+(\xi) \leqslant
\sum_{M\in {\mathfrak S}} \tNN_{M,Q(1+O(Q^{-\alpha/2}))}^+ (\xi)+O (Q^{1+\alpha} \ln Q) ,
\end{equation}
\begin{equation*}
\sum_{M\in {\mathfrak S}} \tNN_{M,Q}^+(\xi) \leqslant
\sum_{M\in {\mathfrak S}} \NN_{M,Q(1+O(Q^{-\alpha/2}))}^+ (\xi)+O (Q^{1+\alpha} \ln Q).
\end{equation*}

{\em (iii)} For $M\in \SSS$ we have individually
\[
\begin{split}
\NN_{M,Q}(\xi) & \leqslant
 \tNN_{M,Q(1+O(Q^{-1}))} (\xi)+O (Q^{1+\alpha}), \\ \tNN_{M,Q}(\xi) & \leqslant
 \NN_{M,Q(1+O(Q^{-1}))} (\xi)+O (Q^{1+\alpha}).
 \end{split}
\]
\end{lemma}

\begin{proof}
(i) Assume first $\frac{|a|}{c}>k>\vert \gamma \omega\vert$ or $\frac{|a|}{c}<k<\vert \gamma \omega\vert$. Since
$\Psi(\gamma)$ is the $x$-intercept of the geodesics from $\omega$ to $\gamma \omega$, it follows that $|\gamma \omega |<k $
if and only if $|\Psi(\gamma)|<k$, so we have (assuming $a>0$, the other case being similar with $k$ replaced by $-k$ below)
\[
\left|\frac{a}{c}-k \right| <
\left| \Psi(\gamma)-\frac{a}{c}\right| <
\left| \Phi(\gamma)-\frac{a}{c}\right| +\big|\Psi(\gamma)-\Phi(\gamma) \big|\ll \frac{1}{c^2} ,
\]
by \eqref{4.2}
and $|\Phi(\gamma)-\frac{a}{c}|=\frac{1}{c^2}\frac{|u+d/c|}{(u+d/c)^2+v^2}\ll\frac{1}{c^2}$. Since $k\in
\Q$, there are finitely many such pairs $(a,c)$, and repeating the argument with $\frac{a}{c}$ replaced by
$\frac{b}{d}$ we obtain that there are finitely many such matrices $\gamma$ (also using the fact that there are finitely
many $\gamma$ with $k$ between $\frac{|a|}{c}$ and $\frac{|b|}{d}$). From \eqref{7.70} and the fact that $\Phi (\gamma)\in \Q$,
there are also finitely many matrices $\gamma'$ satisfying the assumptions.

Finally assume $\frac{|a|}{c}>k>\vert \gamma' \omega\vert$, the remaining case being similar. Then
$|\Phi(\gamma')|<k<\frac{|a|}{c}$, and as before we have (assuming $a>0$, otherwise replace $k$ by $-k$)
\[
\left| \frac{a}{c}-k \right| <
\left| \Phi(\gamma')-\frac{a}{c}\right| < \left| \Phi(\gamma)-\frac{a}{c}\right| +
\vert \Phi(\gamma)-\Phi(\gamma') \vert \ll \frac{1}{c^2} ,
\]
and we conclude as in the previous case.

(ii) Next we look at pairs $(\gamma,\gamma'=\gamma M)$ satisfying \eqref{7.70}, with $\vert a\vert \leqslant k c$ and
$\| \gamma\|,\|\gamma M \| \leqslant Q$, estimating their contribution to
the left-hand side of \eqref{7.5} according to whether $Y_\gamma<Q^{2\alpha}$ or $Y_\gamma\geqslant
Q^{2\alpha}$. By part (i) we can assume $X_\gamma\ll Y_\gamma$, and by \eqref{7.3}
we have $\|\gamma\|\ll Q^{\alpha}$ or $\|\gamma\|\gg Q^{\alpha}$, respectively in the two cases.

Assume first $Y_\gamma<Q^{2\alpha}$. With the Farey interval $J_\gamma$ associated
to $\gamma$ defined before Lemma \ref{L4.1}, the Farey intervals $J_{\gamma}$ and $J_{\gamma M}$ are
either disjoint or one contains the other. Since each Farey interval is associated with at most three matrices
$\gamma\in \Gamma$, it follows as in the proof of \cite[Proposition 3]{BPPZ} that the number of
pairs $(\gamma, \gamma M)$ is $\ll Q^{1+\alpha}\ln Q$.

Therefore we are left to consider pairs $(\gamma,\gamma M)$ with $Y_{\gamma}\geqslant Q^{2\alpha}$,
$\frac{|a|}{c}, |\gamma \omega|,|\gamma M\omega|\le k$, and with
 $\|\gamma\|,\|\gamma M\| \ll Q$. These conditions are satisfied by $\gamma$ in either $\SS_{M,Q}(\xi)$
or $\widetilde{\SS}_{M,Q}(\xi)$, as we can assume $X_\gamma\ll Y_\gamma$, $X_{\gamma M}\ll Y_{\gamma M}$ by part
(i). Without loss of generality we assume $c>d>0$ (otherwise substitute $a$ in terms of $b,c,d$ in the left
side of \eqref{7.3}), and show that
\begin{equation}\label{7.6}
\frac{\vert (2a-uc)(dY_M+cZ_M)-Y_M\vert}{c^2}  \ll \frac{Q^2}{c} \leqslant Q^{2-\alpha} .
\end{equation}
Indeed, the inequality $\| \gamma M\| \ll Q$ plainly gives $Y_M\ll \frac{Q^2}{d^2}$.
Employing also $(dY_M+cZ_M)^2+v^2 c^2=Y_{\gamma M} Y_M$ we arrive at \eqref{7.6}. By \eqref{7.3} the claim
follows.

(iii) By proof of part (ii) it remains to consider pairs $(\gamma,\gamma M)$ with $Y_{\gamma}\geqslant
Q^\alpha$, $|\gamma \omega|,|\gamma M\omega|\le k$, and with $d<0$. Without loss of generality we can assume
$c>|d|$ (otherwise substitute $a$ in terms of $b,c,d$ in the left side of \eqref{7.3}), and \eqref{7.6}
follows trivially since $M$ is fixed, with the upper bound being now $Q^{-\alpha}$ .
\end{proof}

\begin{lemma}\label{L7.3}
For any $M\in {\mathfrak S}$, uniformly in $M$ and $\xi$,
\begin{equation*}
\tNN_{M,Q}(\xi)  =\frac{Q^2}{\zeta(2)} \operatorname{Vol} (S_{M,\xi})+O_\varepsilon (Q^{11/6+\varepsilon}).
\end{equation*}
\end{lemma}

\begin{proof}
$\tNN_{M,Q}(\xi)$ represents the sum over $c\in \{ 1,\ldots,Q\}$ of the number of integer lattice points
$(a,d)$ with $ad\equiv 1 \pmod{c}$ in the region $\Omega =\Omega_{M,Q,c}(\xi)$ of points $(a,d)\in [-kc,kc]\times
[-\widetilde{Q},\widetilde{Q}]$ for which $\vert \Xi_M (c,d)\vert \leqslant \frac{\xi}{Q^2}$ and
$\big( 1+\frac{a^2}{c^2}-2u\frac{a}{c}\big) \max\{ k^2c^2+d^2+2cdu,c^2 X_M +d^2 Y_M +2cdZ_M\} \leqslant Q^2$.
Applying Lemma 7 of \cite{BPPZ} with $q=c$, $A(\Omega ) \leqslant cQ$,
$\ell (\partial \Omega ) \ll Q$, and $L=c^{5/6}$, we find
\begin{equation*}
\tNN_{M,Q}(\xi) =\sum_{c=1}^Q \left( \frac{\varphi (c)}{c^2} A \big(
\Omega_{M,Q,c} (\xi )\big) +O_\varepsilon (Qc^{-1/6+\varepsilon})\right) .
\end{equation*}
The function $h(c)=\frac{1}{c} A(\Omega_{M,Q,c}(\xi))$ with
$\| h\|_\infty \leqslant Q$ is continuous and piecewise $C^1$, with the number of critical
points bounded above by a constant independent of $M,c,\omega,\xi$ (see the discussion
following \cite[Eq. (7.14)]{BPPZ}). Applying
M\" obius summation (as in \cite[Lemma 2.3]{BCZ0}) and the change of variables $(c,u,v)=(Qx, Qxz,Qy)$, we find
\begin{equation*}
\begin{split}
\tNN_{M,Q}(\xi) &= \frac{1}{\zeta(2)} \int_0^Q A\big( \Omega_{M,Q,c} (\xi)\big) \frac{dc}{c}
+O_\varepsilon (Q^{11/6+\varepsilon}) \\  &
= \frac{Q^2}{\zeta(2)} \operatorname{Vol} (S_{M,\xi})
+O_\varepsilon (Q^{11/6+\varepsilon}).
\end{split}
\end{equation*}
\end{proof}

When restricting to the subset $\tSS_{M,Q}^+(\xi)$, the following improved estimate holds.
\begin{lemma}\label{L7.4}
Let $\SSS'\subseteq \SSS$ be any subset.
For each $\beta_0 \in ( \frac{1}{2},1)$ the following estimate holds
(uniformly in $\xi$ on compacts), with $S^+_{M,\xi}=\{ (x,y,z)\in S_{M,\xi}: y>0\}$:
\begin{equation}\label{7.7}
\sum_{M\in {\mathfrak S'}} \tNN^+_{M,Q}(\xi)
=\frac{Q^2}{\zeta(2)} \sum\limits_{\substack{M\in {\mathfrak S'} \\ X_M,Y_M \leqslant Q^{2\beta_0}}}
\hspace{-15pt} \operatorname{Vol} (S_{M,\xi}^+) + O_{\xi,\omega} (Q^{(11+\beta_0)/6}).
\end{equation}
\end{lemma}

\begin{proof}
It is convenient to use the dichotomy from the proof of Lemma \ref{L6.1}:

{\bf Case (A).} $\tan \theta<v$, which is equivalent to $d>(1-u)c$ and yields $X_M \ll Y_M \leqslant Q^{2\beta_0}$.
This will be essential in getting good bounds for the error terms below.
The contribution of this case to $\tNN^+_{M,Q}(\xi)$ is
\[
\NN^{\prime}_{M,Q}(\xi) =\sum_{c\leqslant Q/(v\sqrt{X_M})} \# \big\{ (a,d) \in \Z^2 \cap
\Omega^\prime_{M,Q,c,\omega}(\xi),\
ad\equiv 1 \hspace{-5pt} \pmod{c} \big\} ,
\]
where $\Omega=\Omega^\prime_{M,Q,c,\omega}(\xi)$ is the set of points $(a,d)\in [-kc,kc] \times [(1-u)c,\widetilde{Q}]$
with
\begin{equation}\label{7.8}
\begin{cases}
\vert \Xi_M (c,d)\vert \leqslant \frac{\xi}{Q^2} , \\
\max\{k^2c^2+d^2+2cdu, c^2 X_M +d^2 Y_M +2cdZ_M\} \leqslant \frac{c^2}{(a-uc)^2+v^2c^2} Q^2 .
\end{cases}
\end{equation}
The inequalities $(1-u)c \leqslant d \leqslant \frac{Q}{v\sqrt{Y_M}}$ show that if $\Omega \neq \emptyset$ then
$\vert a\vert \leqslant kc\ll_\omega \frac{Q}{\sqrt{Y_M}}$, and so $A(\Omega) \ll_\omega \frac{Qc}{\sqrt{Y_M}}$
and $\ell (\partial \Omega) \ll_\omega c+\frac{Q}{\sqrt{Y_M}} \ll_\omega \frac{Q}{\sqrt{Y_M}} $.
Applying Lemma 7 of \cite{BPPZ} with $q=c$ and $L=c^{5/6}$, we find
\begin{equation}\label{7.9}
\NN^\prime_{M,Q}(\xi) =\sum_{c\leqslant Q/ ((1-u)v\sqrt{Y_M})} \bigg( \frac{\varphi (c)}{c^2} A
\big( \Omega_{M,Q,c,\omega}^{\prime} (\xi)\big) +O_\varepsilon \Big( \frac{Q c^{-1/6+\varepsilon}}{\sqrt{Y_M}}
\Big)\bigg) .
\end{equation}
According to Lemma \ref{L6.1} we
should sum over $M$ with $\max\{ X_M,Y_M\} \leqslant Q^{2\beta_0}$.
Since $X_M \ll Y_M$, once the entries $C$ and $D$ of $M$ are fixed, the entries $A$ and $B$ can only take $O(1)$ values.
This helps us to conclude that the total contribution of the error term
in \eqref{7.9} to \eqref{7.7} is
\begin{equation*}
\begin{split}
{\mathcal E}_1\ll_\omega \sum\limits_{\substack{M\in {\mathfrak S} \\ X_M \ll Y_M \leqslant Q^{2\beta_0}}} & \sum_{c\ll Q/\sqrt{Y_M}}  \frac{Qc^{-1/6+\varepsilon}}{\sqrt{Y_M}}
\ll_\omega \sum\limits_{\substack{M\in {\mathfrak S} \\ X_M \ll Y_M \leqslant Q^{2\beta_0}}} \frac{Q}{\sqrt{Y_M}}
\Big( \frac{Q}{\sqrt{Y_M}}\Big)^{5/6+\varepsilon} \\
& \ll Q^{11/6+\varepsilon} \sum_{C^2+D^2 \leqslant Q^{2\beta_0}} (C^2+D^2)^{-11/6} \ll Q^{(11+\beta_0)/6+\varepsilon} .
\end{split}
\end{equation*}
The function $h(c)=\frac{1}{c} A (\Omega_{M,Q,c,\omega}^\prime (\xi))$ with
$\| h\|_\infty \ll_\omega \frac{Q}{\sqrt{C^2+D^2}}$ is continuous and piecewise $C^1$, with the number of its critical points
bounded by a universal constant independent of $M,c,\omega,\xi$.
M\" obius summation over $c$ applied to $h$ (as in \cite[Lemma 2.3]{BCZ0}), the change of variables
$(c,a,d)=(Qx,Qxz,Qy)$ and Lemma \ref{L6.1} provide, with
\[
{\mathcal E}_2 =\sum\limits_{\substack{M\in {\mathfrak S} \\ X_M \ll Y_M \leqslant Q^{2\beta_0}}}
\frac{Q\ln Q}{\sqrt{Y_M}} \ll_\varepsilon Q^{1+\beta_0+\varepsilon} ,
\]
the estimates
\begin{equation}\label{7.10}
\begin{split}
\sum_{M\in {\mathfrak S}^\prime} \NN^\prime_{M,Q}(\xi) & =
\sum\limits_{\substack{M\in \mathfrak{S}^\prime \\ X_M,Y_M \leqslant Q^{2\beta_0}}}
\sum\limits_{c\leqslant (1-u)Q/(v\sqrt{Y_M})} \frac{\varphi (c)}{c} h(c)
+ {\mathcal E}_1 \\
&  = \frac{1}{\zeta(2)}
\sum\limits_{\substack{M\in {\mathfrak S}^\prime \\ X_M,Y_M \leqslant Q^{2\beta_0}}}
\int_0^{\frac{(1-u)Q}{v\sqrt{Y_M}}}  A \big( \Omega_{M,Q,c,\omega}^{\prime}
(\xi)\big)\frac{dc}{c} + {\mathcal E}_2 +{\mathcal E}_1 \\
& =\frac{Q^2}{\zeta(2)}  \sum\limits_{\substack{M\in {\mathfrak S}^\prime \\ X_M,Y_M \leqslant Q^{2\beta_0}}}
\operatorname{Vol}(S^\prime_{M,\xi}) +O_{\varepsilon,\omega,\xi} (Q^{(11+\beta_0)/6+\varepsilon}) ,
\end{split}
\end{equation}
where $S^\prime_{M,\xi}$ denotes the subset of $S^+_{M,\xi}$ with the additional condition $(1-u)x\leqslant y$.

{\bf Case (B).} $\tan \theta >v$, which is equivalent to $d<(1-u)c$ and yields $Y_M \ll X_M \leqslant Q^{2\beta_0}$.
Now the contribution $\NN^{\prime\prime}_{M,Q}(\xi)$ of this case to $\widetilde{\NN}^+_{M,Q}(\xi)$ is obtained by counting
integer triples $(a,c,d)$ with $ad\equiv 1  \pmod{c}$ in the region
$\Omega=\Omega^{\prime\prime}_{M,Q,c,\omega}(\xi)$ defined by $(a,d) \in [-kc,kc]\times (0,(1-u)c]$ and by \eqref{7.8},
and which satisfies $A(\Omega) \ll_\omega c^2$ and $\ell (\partial \Omega) \ll_\omega c$ with $c\leqslant \frac{Q}{\sqrt{X_M}}$.
Proceeding as in Case 1 but with the roles of $X_M$ and $Y_M$ reversed, we find
\begin{equation*}
\sum_{M\in {\mathfrak S}^\prime} \NN^{\prime\prime}_{M,Q}(\xi) =
\frac{Q^2}{\zeta(2)}  \sum\limits_{\substack{M\in {\mathfrak S}^\prime \\ X_M,Y_M \leqslant Q^{2\beta_0}}}
\operatorname{Vol}(S^{\prime\prime}_{M,\xi}) +O_{\varepsilon,\omega,\xi} (Q^{(11+\beta_0)/6+\varepsilon}) ,
\end{equation*}
with $S_{M,\xi}^{\prime\prime}$ denoting the subset of $S_{M,\xi}^+$ with $y\leqslant (1-u)x$.
\end{proof}

\begin{corollary}\label{C7.5}
For each $\beta_0 \in ( \frac{1}{2},1)$ the following estimate holds:
\begin{equation*}
\sum_{M\in {\mathfrak S}} \tNN^+_{M,Q}(\xi)
=\frac{Q^2}{\zeta(2)} \sum\limits_{M\in {\mathfrak S}}
 \operatorname{Vol} (S_{M,\xi}^+) + O_{\xi,\omega} (Q^{(11+\beta_0)/6}).
\end{equation*}
\end{corollary}
\begin{proof}
This follows from an adaptation of (7.17) in \cite{BPPZ}, Lemma \ref{L7.4}, and
\[
\sum_{C^2+D^2 \geqslant Q^\sigma} (C^2+D^2)^{-2} \ll \int_{Q^{\sigma/2}}^\infty \frac{dr}{r^3}
=\textstyle\frac{1}{2} Q^{-\sigma}, \qquad \sigma >0 .
\]

For $S^\prime_{M,\xi}$ we are
in Case (A). Therefore $y^2 \ll \frac{1}{Y_M} \ll \frac{1}{T}$ and $x<\frac{y}{1-u}$,
giving $r^2=x^2+y^2 \ll \frac{1}{T}$. On the other hand the proof of Lemma \ref{L6.1}
(there is no need to assume the integrality of $c$ and $d$ for this estimate) provides
$\frac{x}{ux+y}=\omega_M +O(\frac{1}{T})$ with $\omega_M =\frac{Z_M-uY_M}{v^2 Y_M}$.
Projecting $S_{M,\xi}^\prime$ on the $(x,y)$-coordinates we find
\[
\operatorname{Vol}(S_{M,\xi}^\prime) \leqslant 2k A\Big(
\Big\{ (x,y)\in (0,\infty)^2: x^2+y^2 \ll_\omega T^{-1}, \frac{x}{ux+y} =\omega_M +O(T^{-1}) \Big\}\Big).
\]
Changing coordinates to $x=\alpha$, $ux+y=\beta$, the region in the right-hand side above
is mapped onto the intersection of the ellipse $\alpha^2 +(\beta-u\alpha)^2 \ll T^{-1}$ with
$\alpha^2+\beta^2 \ll_\omega T^{-1}$ and the wedge $\frac{\alpha}{\beta} =\omega_M+O(\frac{1}{T})$,
$\alpha,\beta >0$. Using polar coordinates it is immediate that its area is $\ll_\omega T_M^{-2}=T^{-2}\ll \| M\|^{-4}$.

The situation of $S^{\prime\prime}_{M,\xi}$ is similar.
\end{proof}

\section{Extra symmetries for $\omega=i$}\label{sec8}

In this section we assume $\omega=i$ and make use of extra symmetries of the hyperbolic lattice centered at $i$.
For each matrix $M=\left( \begin{smallmatrix} A & B \\ C & D \end{smallmatrix}\right)
\in\Gamma$ consider the cardinality $\NN_{M,Q}(\xi)$ of the set $\SS_{M,Q}(\xi)$ defined in \eqref{1.6}. We first
show that in the expression \eqref{1.5} we can restrict to matrices $M$ having positive entries.
Then we show that, except for a set of cardinality $\ll_\xi 1$ of such matrices $M$, we can restrict to counting only
elements $\gamma\in \SS_{M,Q}(\xi)$ with positive elements in the second row. Thus we can make use of
results of the previous sections to estimate the quantity $\RR^\Phi_Q (\xi)$.

Recall the set ${\mathfrak S}$ of matrices $M\neq I$ with nonnegative entries, and denote $s=\left(
\begin{smallmatrix} 0 & -1 \\ 1 & 0 \end{smallmatrix}\right)$. The signs of the entries of $sM$, $Ms$, and $sMs$
show that $\Gamma \setminus \{ I,s\}$ is partitioned into ${\mathfrak S}\cup s{\mathfrak S} \cup {\mathfrak S}s
\cup s{\mathfrak S}s$. The equalities $\Phi (gs)=\Phi (g)$, $si=i$, and $\| g\|=\| gs\|$ yield
\begin{equation*}
\begin{split}
& \SS_{Ms,Q}(\xi)=\SS_{M,Q}(\xi),\quad \SS_{sM,Q}(\xi)=\SS_{M,Q}(\xi) s,
\\  &
\NN_{M,Q}(\xi)=\NN_{Ms,Q}(\xi)=\NN_{sM,Q}(\xi)=\NN_{sMs,Q}(\xi).
\end{split}
\end{equation*}
Therefore \eqref{1.5} becomes
\begin{equation*}
\RR^\Phi_Q (\xi) =2 \sum_{M\in {\mathfrak S}} \NN_{M,Q}(\xi).
\end{equation*}

Let $\SS^+_{M,Q}(\xi)$, respectively $\SS^-_{M,Q}(\xi)$, denote the subsets of $\SS_{M,Q}(\xi)$ consisting of
matrices with $c,d>0$, respectively $c>0>d$, of cardinality $\NN^\pm_{M,Q}(\xi)$ respectively.

Recall the finite set $\FF(\xi)$ defined in \eqref{5.1} and let
\begin{equation*}
\widetilde{\FF}(\xi)  : =\FF (\xi) \cup s\FF(\xi) \cup s\FF(\xi)^{-1} .
\end{equation*}

\begin{lemma}\label{L8.1}
$\mathrm{(i)}$ For any $M\in {\mathfrak S}$ the mapping $\gamma \mapsto
\gamma M s^{-1}$ defines a bijection between the sets
$$\SS^{-, \ast}_{M,Q}(\xi)=\left\{ \gamma=\Big( \begin{smallmatrix} a & b \\ c & d
\end{smallmatrix}\right)\in \SS^-_{M,Q}(\xi): \textstyle\frac{B}{D} < -\frac{d}{c} < \frac{A}{C}\Big\} \quad
\mbox{\rm and} \quad \SS^+_{M^t,Q}(\xi).$$

$\mathrm{(ii)}$ If $M\in {\mathfrak S} \setminus \widetilde{\FF} (\xi)$, then
$\SS_{M,Q}^{-}(\xi)=\SS_{M,Q}^{-,\ast}(\xi)$, so that $\NN^-_{M,Q}(\xi)=\NN^+_{M^t,Q}(\xi)$.
\end{lemma}

\begin{proof}
$\mathrm{(i)}$ Let $\gamma= \left( \begin{smallmatrix} a & b \\ c & d
\end{smallmatrix}\right) \in \SS_{M,Q}^{-,\ast}(\xi)$ and $\gamma M s^{-1}=\left( \begin{smallmatrix} * & * \\
\tilde{c} & \tilde{d} \end{smallmatrix}\right)$. On one hand $c>0$ and $\frac{B}{D}<\frac{-d}{c} <
\frac{A}{C}$ imply $\tilde{c}=-cB-dD >0, \tilde{d}=cA+dC >0$, while $\tilde{c}>0$ and $\tilde{d}>0$ imply
$c=\tilde{c}C+\tilde{d}D>0$ and $d=-\tilde{c}A-\tilde{d}B <0$. On the other hand, utilizing also
$sM^{-1} s^{-1} =M^t$, $\Phi(gs^{-1})=\Phi(g)$, $\|  \gamma M s^{-1}M^t \|=\| \gamma\|=\| \gamma s^{-1}\|$, we conclude
that the map above is a bijection.

$\mathrm{(ii)}$ Suppose, by contradiction, that there exists $\gamma\in \SS_{M,Q}^{-}(\xi) \setminus
\SS_{M,Q}^{-,\ast}(\xi)$. Setting
$$\gamma M =\left( \begin{matrix} aA+bC & aB+bD \\ cA+dC & cB+dD
\end{matrix}\right)=\left( \begin{matrix} a^\prime & b^\prime \\ c^\prime & d^\prime
\end{matrix}\right)=\gamma^\prime,$$
notice that $Q\geqslant \max\{ c,-d, \vert c^\prime\vert,\vert d^\prime\vert\}$ and two cases can occur:

1) $\  \frac{-d}{c} <\frac{B}{D} <\frac{A}{C}$. Then $c^\prime >0$, $d^\prime >0$, and
$\frac{b^\prime}{d^\prime} < \frac{a^\prime}{c^\prime} \leqslant \frac{a}{c} <
\frac{-b}{-d}$. By Lemma \ref{L5.2}, $\gamma s^{-1}=\gamma^\prime M_0$ with $M_0\in\FF(\xi)$, so $M =(M_0 s)^{-1}$,
contradiction.

2) $\  \frac{B}{D}<\frac{A}{C} <\frac{-d}{c}$. Then $c^\prime <0$, $d^\prime <0$, and
$\frac{a}{c} <\frac{-b}{-d} \leqslant \frac{-b^\prime}{-d^\prime} < \frac{-a^\prime}{-c^\prime}$.
Lemma \ref{L5.2} gives $-\gamma^\prime=\gamma s^{-1} M_0$ with $M_0\in \FF(\xi)$, so $M=-s^{-1} M_0$, contradiction.
\end{proof}

\begin{remark}\label{R8.2} The analogue of Lemma \ref{L8.1} (ii) holds for $\tSS_{M,Q}^{-}(\xi)$ in place of
$\SS_{M,Q}^-(\xi)$ (see the notation preceding Lemma \ref{L7.2}). Namely, there is no matrix $\gamma=\left(
\begin{smallmatrix} a & b \\ c & d
\end{smallmatrix}\right)\in \tSS_{M,Q}^{-}(\xi)$ with either $-\frac{d}{c}
<\frac{B}{D}$ or $\frac{A}{C} <\frac{-d}{c}$. Indeed, referring to the notation in the proof of Lemma \ref{L8.1}
(ii), we have $c^{\prime 2}+d^{\prime 2} =c^2 X_M +d^2 Y_M +2cd Z_M\le Q^2$, so that  $Q\geqslant \max\{ c,-d,
\vert c^\prime\vert,\vert d^\prime\vert\}$, and the rest of the proof goes through unchanged.
\end{remark}

Consider now the region $S_{M,\xi}$ in \eqref{1.8}, which in the present case $\omega=i$
becomes the region of triples $(x,y,z)\in [0,1]\times [-1,1]^2$ for which\footnote{As
$X_M,Y_M\geqslant 1$, $Z_M \geqslant 0$, when $y>0$ the inequality $x^2+y^2 \leqslant \frac{1}{1+z^2}$ is
obsolete.}
\begin{equation}\label{8.1}
 \vert \Xi_M (x,y)\vert \leqslant \xi,\quad
\max\{ x^2+y^2 , x^2 X_M +y^2 Y_M +2xy Z_M\} \leqslant \frac{1}{1+z^2} .
\end{equation}
Consider also the subsets $S_{M,\xi}^+$, $S_{M,\xi}^-$, $S_{M,\xi}^{-,\ast}$ of $S_{M,\xi}$ defined, respectively, by
$y>0$, $y<0$, $y<0$ and $\frac{B}{D}< \frac{-y}{x}<\frac{A}{C}$. As in Lemma \ref{L8.1},
the mapping $(x,y) \mapsto (x,y)Ms^{-1}=(\tilde{x},\tilde{y})$ defines a diffeomorphism between the
sets
$S_{M,\xi}^{-,\ast}$ and $S_{M^t,\xi}^+$, showing in particular that
\begin{equation}\label{8.2}
\operatorname{Vol} (S_{M,\xi}^{-,\ast}) = \operatorname{Vol} (S_{M^t,\xi}^+),\qquad \forall M\in {\mathfrak S}
\setminus \widetilde{\FF}(\xi) .
\end{equation}

\begin{lemma}\label{L8.3}
If $M\in {\mathfrak S} \setminus \widetilde{\FF}(\xi)$, then $\operatorname{Vol}(S^-_{M,\xi}) =
\operatorname{Vol} (S^{-,\ast}_{M,\xi} )$.
\end{lemma}

\begin{proof}
Suppose,  by contradiction, that the set $S_{M,\xi}^- \setminus S_{M,\xi}^{-,\ast}$ contains an interior point
$(x_0,y_0,z_0)$, that is
$(x_0,y_0,z_0)\in (0,1)\times (-1,0)\times (-1,1)$ satisfies both \eqref{8.1} and
$\frac{-y}{x} \in ( 0,\frac{B}{D}) \cup ( \frac{A}{C},\infty)$. Now the set of rational points
\[
\Omega_Q=\Big\{\Big(\frac{c}{Q}, \frac{d}{Q}, \frac{a}{Q} \Big): (c,d)=1,\  Q\geqslant c>0>d\geqslant -Q,\
ad\equiv 1 \hspace{-5pt}\pmod{c},\  |a|<c\Big\}
\]
is dense in $D=[0,1]\times [-1,0]\times [-1,1]$; indeed for each parallelepiped $R \subseteq D$, we can
count the number of points in the scaled sets $QR\cap Q\Omega_Q$, as in the proof of Lemma \ref{L7.1}, and conclude
that
$\Omega_Q\cap R$ is dense in $R$. Therefore for large enough $Q$, we can find points in $\Omega_Q$ arbitrarily
close to
$(x_0,y_0,x_0z_0)$, and it follows that there exists $(a,c,d) \in \tSS_{M,Q} (\xi)$ with
$c>0>d$ and $\frac{-d}{c} \in ( 0,\frac{B}{D}) \cup ( \frac{A}{C},\infty)$,
which contradicts Remark \ref{R8.2}.
\end{proof}

Estimates for $\NN^-_{M,Q}(\xi)$ with $M\notin \widetilde{\FF} (\xi)$
are derived from those on $\NN^+_{M^t,Q} (\xi)$ by Lemma \ref{L8.1}. We can now estimate $\sum_{M\in {\mathfrak S}}
\NN_{M,Q} (\xi)$ by first breaking the sum into sums over $\widetilde{\FF}(\xi)$ and over ${\mathfrak
S}\backslash\widetilde{\FF}(\xi)$; for the first sum we use Lemma \ref{L7.3}, while for the second we use
$\NN^-_{M,Q}(\xi)=\NN^+_{M^t,Q} (\xi)$ and Lemmas \ref{L7.2} and \ref{L7.4}. Finally, employing \eqref{8.2} and Lemma
\ref{L8.3}, we find
\begin{equation}\label{8.3}
\sum_{M\in {\mathfrak S}} \NN_{M,Q} (\xi) =\frac{Q^2}{\zeta(2)}
\sum_{M\in {\mathfrak S}} \operatorname{Vol} (S_{M,\xi}) + O_{\xi} (Q^{(11+\beta_0)/6}).
\end{equation}
To complete the sum to $M\in \Gamma$, note that $\operatorname{Vol} (S_{M,\xi})=\operatorname{Vol} (S_{Ms,\xi})$.
This is also seen to coincide with $\operatorname{Vol} (S_{sM,\xi})$, and thus with
$\operatorname{Vol} (S_{sMs,\xi})$, by employing $\Xi_{sM} (y,-x)=\Xi_M (x,y)$ and the change of variable
$(x,y)\mapsto (-y,x)$ if $(x,y)\in [0,1]\times [-1,0]$, respectively
$(x,y) \mapsto (y,-x)\mapsto (y,-x)$ if $(x,y)\in [0,1]\times [0,1]$. This proves \eqref{1.9}.

\section{A closed form formula for $\operatorname{Vol}(S_{M,\xi})$}\label{sec9}

In this section we evaluate the volume of the body $S_{M,\xi}$ in \eqref{1.8} for \emph{arbitrary} $\omega$,
which leads to the formula in Conjecture \ref{Conj1}. For $\omega=i$, the proof of this conjecture is based on the results of
the previous section.

The volume can be brought in closed form using the substitution
\begin{equation}\label{9.1}
x\omega+y=re^{i\theta},\quad z=v\tan t+u ,
\end{equation}
with $\Xi_M(x,y)$ given by the first Eq. \eqref{2.5}, $k^2 x^2+y^2+2uxy=r^2$, $x^2 X_M
+y^2 Y_M +2xy Z_M =r^2 \sinh \ell (M)
( \coth \ell(M)+\cos (\theta_{M} -2\theta))$, $k^2+z^2-2uz=\frac{v^2}{\cos^2 t}$, to
\begin{equation}\label{9.2}
\operatorname{Vol} (S_{M,\xi})= v\int_{\arctan ((-k-u)/v)}^{\arctan ((k-u)/v)} B_M
(\xi,t)\frac{dt}{\cos^2 t}
= v\int_{\beta/2-\pi/2}^{\beta/2} B_M (\xi,t)\frac{dt}{\cos^2 t} ,
\end{equation}
with $\beta\in(0,\pi)$ such that $\omega=ke^{i\beta}$, where $B_M (\xi,t)$ is the area of
the region defined in polar coordinates $(r,\theta)$ by
\begin{equation}\label{9.3}
\begin{cases}
\vspace{0.2cm}
\frac{r}{v} (\sin\theta, k\sin(\beta-\theta) ) \in \big[ 0,\frac{\wQ}{kQ}\big] \times \big[
-\frac{\wQ}{Q},\frac{\wQ}{Q} \big] \\
\frac{v}{\xi} \frac{\vert \sin (\theta_M -2\theta)\vert}{U_M+\cos (\theta_M -2\theta)}
\leqslant r^2 \leqslant \frac{\cos^2 t}{v^2} \min\Big\{ 1,\frac{1}{\sinh \ell(M) (U_M+
\cos (\theta_{M}-2\theta))} \Big\} ,
\end{cases}
\end{equation}
with $U_M=\coth \ell (M)=\frac{T}{\sqrt{T^2-\Delta^2}} >1$, $T=\| M\|^2$.

Using the second condition in \eqref{9.3}, we have $r^2 \leqslant \frac{1}{v^2} \leqslant
\frac{1}{k(k-|u|)}$. Hence the first condition in \eqref{9.3} can be replaced by
$0\leqslant \theta \leqslant \pi$, and the area $B_M(\xi,t)$ can be expressed in closed form,
with $f_+=\max \{ f, 0\}$, as
\begin{equation*}
\frac{1}{2v}\int_{0}^{\pi}
\frac{\big( \frac{\cos^2 t}{v^2} \min\big\{ \frac{1}{\sinh \ell (M)},U_M +\cos (\theta_M -2\theta)\big\}
-\frac{v}{\xi} \vert \sin (\theta_M -2\theta)\vert\big)_+}{U_M+\cos (\theta_M -2\theta)}\
d\theta.
\end{equation*}

Since we are interested in the pair correlation of the angles $\theta_\gamma$, we
define
\[
\RR_Q^\theta(\xi):=\#\Big\{(\gamma,\gamma^\prime)\in \RRR_Q^2 :
\gamma\ne\gamma^\prime,
\ 0\leqslant \theta(\gamma)-\theta(\gamma')<\frac{\xi}{Q^2} \Big\}.
\]

Following the approximation arguments from Section 8 of \cite{BPPZ}, from \eqref{8.3} we obtain the following
asymptotics:
\begin{proposition}\label{P9.1}
For $\omega=i$ one has
 \begin{equation*}
 \begin{split}
\RR_Q^\theta (\xi) & =\frac{Q^2}{2\zeta(2)} \sum_{M\in\Gamma\backslash\{ I\}} B_M(\xi)
+O_{\xi,\varepsilon} (Q^{47/24+\varepsilon}),\quad \mbox{where} \\
B_M (\xi) & =v\int_{\beta/2-\pi/2}^{\beta/2} B_M \Big( \frac{v\xi}{2 \cos^2 t},t \Big)
\frac{dt}{\cos^2 t}.
\end{split}
\end{equation*}
\end{proposition}
Since $B_M ( \frac{\xi}{\cos^2 t},t) = B_M (\xi,0) \cos^2 t$, we find
\[B_M(\xi)= B_M \Big(\frac{v\xi}{2},0\Big) \frac{\pi v}{2} .
\]

Taking derivatives we obtain
\begin{equation*}
B_M^\prime (\xi) =\frac{\pi}{2\xi^2} \int_{I_{\xi,M}}
\frac{\vert \sin (\theta_M -2\theta)\vert}{U_M+\cos (\theta_M -2\theta )}  d\theta ,
\end{equation*}
with $I_{\xi,M} =\{ \theta \in [0,\pi] :
\vert \sin (\theta_M -2\theta)\vert \leqslant \frac{\xi}{\Delta}
\min \{ U_M+\cos (\theta_M -2\theta), \frac{1}{\sinh \ell (M)}\} \}$ (recall
$\Delta=2v^2$).
With $C_M:=\sqrt{\frac{T-\Delta}{T+\Delta}}
=\tanh ( \frac{\ell(M)}{2}) \in (0,1)$ we have
$\sinh \ell (M)=\frac{2C_M}{1-C_M^2}$, $\cosh \ell (M)=\frac{1+C_M^2}{1-C_M^2}$,
$U_M-C_M=\frac{1}{\sinh\ell (M)}$, $\sinh ( \frac{\ell(M)}{2})=\frac{\sqrt{T-2v^2}}{2v}
=\frac{C_M}{\sqrt{1-C_M^2}}$.
Using the change of variable $u=2\theta-\theta_M\in [-\pi,\pi]$
the integrand is even on $[-\pi,\pi]$, and so we have
\begin{equation*}
B_M^\prime (\Delta\xi) =\frac{\pi}{ 2\Delta^2\xi^2} \int_{J_{\xi,M}}
\frac{\sin u}{U_M+\cos u} du ,
\end{equation*}
with $J_{\xi,M} = (J^{(1)}_{\xi,M}\cup J^{(2)}_{\xi,M}) \cap [0,\pi]$, where
\[
\begin{split}
& J^{(1)}_{\xi,M}=\Big\{ u: \cos u \geqslant -C_M, \ \sin u\leqslant
\frac{\xi}{ \sinh \ell (M)} \Big\},\\  &
J^{(2)}_{\xi,M}=\big\{ u: \cos u \leqslant -C_M,\
\sin u \leqslant \xi ( U_M+\cos u ) \big\}.
\end{split}
\]
A direct calculation provides
\begin{equation*}
\begin{split}
J^{(1)}_{\xi,M} & = \begin{cases}
[0,\arccos (-C_M)] \qquad \mbox{\rm if $\xi \geqslant \sinh \ell (M)
=\frac{\sqrt{T^2-\Delta^2}}{\Delta},$} \\
\big[ 0,\arcsin ( \frac{\xi}{\sinh \ell (M)}) \big] \cup
\big[ \pi-\arcsin ( \frac{\xi}{\sinh \ell (M)}), \arccos (-C_M)\big] \\
\hspace{5cm}
\mbox{\rm if $2 \sinh (\frac{\ell(M)}{2}) \leqslant \xi \leqslant
\sinh \ell(M),$} \\
\big[ 0,\arcsin ( \frac{\xi}{\sinh \ell (M)})\big] \qquad
\mbox{\rm if $\xi \leqslant 2 \sinh(
\frac{\ell(M)}{2})=\frac{\sqrt{T-\Delta}}{v}.$}
\end{cases} \\
J^{(2)}_{\xi,M} & = \begin{cases}
[\arccos (-C_M),\pi] \qquad \mbox{\rm if  $\xi\geqslant \sinh \ell(M),$} \\
[\arccos (-C_M),\alpha+\arcsin (U_M \sin\alpha)] \cup
[\pi +\alpha -\arcsin (U_M \sin\alpha),\pi] \\
\hspace{4cm} \mbox{\rm if $2\sinh ( \frac{\ell(M)}{2})
\leqslant \xi \leqslant \sinh \ell(M),$} \\
[\pi+\alpha -\arcsin (U_M \sin\alpha),\pi] \qquad
\mbox{\rm if $\xi \leqslant 2 \sinh ( \frac{\ell(M)}{2}) .$}\end{cases} \\
J_{\xi,M} & = \begin{cases}
[0,\pi] \qquad \mbox{\rm if  $\xi\geqslant \sinh \ell(M),$} \\
\big[ 0,\arcsin(\frac{\xi}{\sinh\ell(M)})\big] \cup
\big[ \pi-\arcsin( \frac{\xi}{ \sinh\ell(M)}),\alpha+\arcsin (U_M \sin\alpha)\big]
\\ \quad \cup\
[\pi +\alpha -\arcsin (U_M \sin\alpha),\pi]
\qquad \mbox{\rm if $2\sinh ( \frac{\ell(M)}{2})
\leqslant \xi \leqslant \sinh \ell(M),$} \\
\big[ 0,\arcsin(\frac{\xi}{\sinh\ell(M)})\big] \cup
[\pi+\alpha -\arcsin (U_M \sin\alpha),\pi] \
\mbox{\rm if $\xi \leqslant 2 \sinh ( \frac{\ell(M)}{2}) $},
\end{cases}
\end{split}
\end{equation*}
where $\alpha=\alpha(\xi)=\arcsin( \frac{\xi}{\sqrt{\xi^2+1}}) \in ( 0, \frac{\pi}{2})$.
With $f_\xi (\ell)$ as in \eqref{1.2} we obtain
$$
B_M'(\Delta \xi)=\frac{\pi}{\Delta^2\xi^2}f_\xi \big( \ell(M)\big).
$$
Letting $R_2^\theta(\xi)=\lim_{Q\rightarrow \infty} \frac{1}{Q^2} \RR_Q^{\theta}(\xi)$,
from Proposition \ref{P9.1} we infer
\[
\frac{dR_2^\theta}{d\xi} (\Delta\xi)=\frac{\pi}{2\zeta(2)\Delta^2\xi^2}\sum_{M\in\Gamma}
f_\xi \big( \ell(M)\big)
\]
(note that $f_\xi (0)=0$ so we can include $I$ in the range of summation).
Taking into account that the pair correlation distribution $R_2(\xi)$ in the introduction
involves normalized angles, and that $B_Q\sim \frac{3}{\Delta}Q^2$, we have
\[
R_2\Big(\frac{3\xi}{\pi\Delta}\Big)=\frac{\Delta}{3}R_2^\theta(\xi) , \quad
g_2\Big(\frac{3\xi}{\pi}\Big)=\frac{\pi\Delta^2}{9} \frac{dR_2^\theta}{d\xi}
(\Delta\xi) .
\]
This leads to the formula for $g_2$ stated in Conjecture \ref{Conj1}, and proves Theorem \ref{Thm1.1} when $\omega=i$.

\section{The case $\omega=\rho$}\label{sec10}

In the case of the other elliptic point $\omega=\rho$ one can take advantage of some
other symmetries to prove Conjecture 1. Consider the matrix
$w=\left( \begin{smallmatrix} 1 & -1 \\ 1 & 0 \end{smallmatrix} \right)$ fixing the point $\rho$.
This time we partition the upper half plane $\H$ into three regions, permuted clockwise by $w$:
\begin{equation*}
\begin{split}
{\bf I} & = \left\{ z\in \H: \operatorname{Re} z >\textstyle \frac{1}{2}, \vert  z-1\vert < 1\right\}, \quad
{\bf II} =\left\{ z\in\H: \operatorname{Re} z< \textstyle\frac{1}{2} , \vert z\vert <1\right\} , \\
{\bf III}  & =\{ z\in \H:\ \vert z-1\vert >1, \vert z\vert > 1\}.
\end{split}
\end{equation*}
The condition that $\gamma \rho$ belongs to one of these is easily stated
in terms of $(X,Y,Z)$ given by \eqref{2.1}, using the relations $\operatorname{Re} (\gamma \rho)=\frac{Z}{Y}$,
$|\gamma \rho|^2=\frac{X}{Y}$, $|\gamma \rho-1 |^2=1+\frac{X-2Z}{Y}$:
\[
\begin{split}
\gamma \rho\in {\bf I} \iff & 2Z>X,\ 2Z>Y, \quad \gamma \rho\in {\bf II} \iff Y>X,\ Y>2Z, \\ &
\gamma \rho\in {\bf III} \iff X>Y, X>2Z.
\end{split}
\]
\begin{figure}[ht]
\centering
\includegraphics*[scale=1, bb=0 0 230 75]{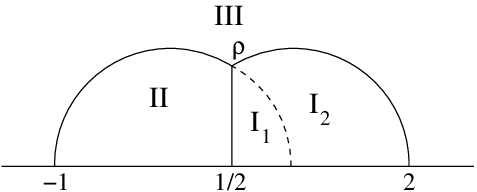}
\caption{Symmetric geodesics through $\rho$}
\label{fig2}
\end{figure}

Next  we determine the restrictions that the condition $\gamma\rho \in {\bf I}$ places on
the entries of $\gamma$. As a consequence of $2Z>Y,2Z>X$, a quick check shows that the
entries of $\gamma$ are
nonzero. Since $abcd=bc+(bc)^2 >0$, we also have that $ac$ and $bd$ have the same sign. In
fact $ac>0$, $bd>0$: if the contrary were true, from $2Z=2(ac+bd)+2ad+1>0$ it would follow
$ad>0$ and without loss of generality we can assume $a>0$, $c<0$, $d>0$, $b<0$; from $ad>-ac-bd$ it would then
follow that $d>-c$, $a>-b$, which implies $ad-bc>1$, a contradiction.
Since $ac>0$, $bd>0$, among the matrices $\gamma, \gamma w, \gamma w^2$ with the same
coordinates $(X,Y,Z)$ precisely one has entries of the same sign. Therefore we assume
from now on that $\gamma$ has positive entries whenever $\gamma \rho \in {\bf I}$.

Using the substitution $b=\frac{ad-1}{c}$, one checks that
\[
2Z>Y \iff 2a-c> \frac{2d+c}{c^2+d^2+cd} \iff 2a>c
\]
(the last equivalence follows since the fraction is less than 1 and $2a-c$ is integral).
Similarly $2Z>X$ if and only if $2d>b$.
In conclusion, the point $\gamma \rho$ belongs to ${\bf I}$ if and only if, after
perhaps replacing $\gamma$ by $\gamma w$ or $\gamma w^2$, we have
\begin{equation}\label{10.1}
a,b,c,d>0, \quad \frac{a}{c}>\frac{1}{2},\quad \frac{b}{d}< 2.
\end{equation}

Since $w$ permutes the regions ${\bf I}$, ${\bf II}$ and ${\bf III}$ and $w\rho=\rho$, the previous discussion
shows that $\Gamma$ can be partitioned as
\begin{equation*}
\begin{split}
\Gamma \setminus \{ I,w,w^2 \} =\bigcup\limits_{r,s\in\{ 0,1,2\}} w^r {\mathfrak M} w^s ,\quad \mbox{\rm where} \\
{\mathfrak M}:=\bigg\{ M=\bigg( \begin{matrix} A & B \\ C & D
\end{matrix}\bigg)\in\Gamma  :
\ C,D >0 ,\ \frac{1}{2}\leqslant \frac{B}{D}<\frac{A}{C}\leqslant 2 \bigg\}.
\end{split}
\end{equation*}

We can now rewrite the sum \eqref{1.5}.
Since $\| g w\|=\|g\|$ and $\Phi (g w)=\Phi(g)$, we have $\SS_{Mw^s,Q}(\xi)=\SS_{M,Q}(\xi)$. Moreover,
$g\mapsto gw^{-r}$ maps $\SS_{M,Q}(\xi)$ bijectively onto  $\SS_{w^rM,Q}(\xi)$, so $\NN_{w^r
Mw^s,Q}(\xi)=\NN_{M,Q}(\xi)$. We infer
\begin{equation*}
\RR_Q^\Phi (\xi)= \frac{9}{2} \sum_{M\in {\mathfrak M}} \NN_{M,Q}(\xi) .
\end{equation*}

To state the equivalent of Lemma \ref{L8.1}, we further divide region ${\bf  I}$ in two regions ${\bf I_1}$ and ${\bf
I_2}$, according as $|z|<1$ or $|z|>1$ (see Fig. \ref{fig2}). Lemma \ref{L3.1} shows that $\gamma \rho\mapsto
\widetilde{\gamma} \rho$ is a bijection of ${\bf I_1}$ onto ${\bf I_2}$. Let also $\MMM_1$, respectively $\MMM_2$,
denote the subset of $M\in \MMM$ with $M\rho \in {\bf I_1}$, respectively  $M\rho \in {\bf I_2}$.

For  $M=\left( \begin{smallmatrix} A & B \\ C & D \end{smallmatrix}\right)\in \MMM$, let $\SS^{-, 1}_{M,Q}(\xi)$,
respectively $\SS^{-, 2}_{M,Q}(\xi)$ be the sets of those $\gamma=\left( \begin{smallmatrix} a & b \\ c & d
\end{smallmatrix}\right)\in \SS^{-}_{M,Q}(\xi) $ for which  $\frac{B}{D} < -\frac{d}{c} < \frac{A+B}{C+D}$,
respectively $\frac{A+B}{C+D} < -\frac{d}{c} < \frac{A}{C}$. With $\FF(\xi)$ as in in \eqref{5.1} define
$$
\widetilde{\FF} (\xi)= \FF(\xi)\cup w\FF(\xi)\cup w^{-1}\FF(\xi)\cup w\FF(\xi)^{-1}\cup  w^{-1}\FF(\xi)^{-1}.
$$

\begin{lemma}\label{L10.1}
{\em (i)} The map $\gamma \mapsto \gamma M w^{-1}$ is a bijection between $\SS^{-,
1}_{M,Q}(\xi)$ and $\SS^{+}_{wM^{-1}w^2,Q}(\xi)$, and the map $\gamma \mapsto \gamma M w$ is a bijection between
$\SS^{-, 2}_{M,Q}(\xi)$ and $\SS^{+}_{w^2 M^{-1}w,Q}(\xi)$.

{\em (ii)} If $M\in {\mathfrak M} \setminus \widetilde{\FF} (\xi)$, then
$\SS_{M,Q}^{-}(\xi)=\SS_{M,Q}^{-,1}(\xi)\cup \SS_{M,Q}^{-,2}(\xi) $, so that $$\NN_{M,Q}^-=\NN_{wM^{-1}w^2,Q}^{+}
+\NN_{w^2 M^{-1}w,Q}^{+}. $$
\end{lemma}
The proof is very similar to that of Lemma \ref{L8.1} and we leave it as an exercise for the reader.

Let now $\MMM_i(\xi)=\MMM_i\backslash \widetilde{\FF}(\xi)$, $i\in \{1,2\}$, and define the sets
\[
\begin{split}
& \SSS_1(\xi)=\{wM^{-1}w^2, w^2M^{-1}w^2 : M\in  \MMM_1(\xi)\}, \\
& \SSS_2(\xi)=\{w^2M^{-1}w, wM^{-1}w : M\in \MMM_2(\xi)\} ,
\end{split}
\]
both easily checked to be contained in $\SSS$ (namely they contain matrices with positive entries). From
Lemma \ref{L10.1} it follows that
\begin{equation*}
\sum\limits_{M\in {\mathfrak M}} \NN_{M,Q}(\xi)  = \sum_{M\in \MMM \cap \widetilde{\FF}(\xi)} \hspace{-10pt} \NN_{M,Q}(\xi) +
\sum_{M\in\MMM_1(\xi)\cup \MMM_2(\xi)} \hspace{-10pt} \NN_{M,Q}^+(\xi)+\sum_{M\in\SSS_1(\xi)\cup \SSS_2(\xi)} \hspace{-10pt} \NN_{M,Q}^+(\xi).
\end{equation*}

Since now we only sum over matrices with positive entries, the approximation arguments employed in the case
$\omega=i$ also apply here, with $S_{M,\xi}$ defined in \eqref{1.8} and $S^\pm_{M,\xi}$ the subset of
$S_{M,\xi}$ defined by the additional condition $y>0$ or $y<0$, leading to
\[\frac{\zeta (2)}{Q^2} \sum\limits_{M\in {\mathfrak M}} \NN_{M,Q}(\xi) \sim \hspace{-10pt}
\sum_{M\in \MMM \cap \widetilde{\FF}(\xi)} \hspace{-10pt} \Vol (S_{M,\xi}) +
\hspace{-5pt} \sum_{M\in\MMM_1(\xi)\cup \MMM_2(\xi)} \hspace{-15pt} \Vol (S_{M,\xi}^+)+
\hspace{-5pt} \sum_{M\in\SSS_1(\xi)\cup \SSS_2(\xi)} \hspace{-15pt} \Vol (S_{M,\xi}^+).
\]
Using equalities analogous to those of Lemma \ref{L10.1}, for volumes instead of the number of lattice points, we find
\begin{equation}\label{10.2}
\RR_Q^\Phi (\xi) \sim \frac{9Q^2}{2\zeta(2)} \sum_{M\in {\mathfrak M}} \operatorname{Vol} (S_{M,\xi}).
\end{equation}

Finally the sum of volumes in \eqref{10.2} can be extended from ${\mathfrak M}$ to $\Gamma$ since
$\operatorname{Vol}(S_{w^r Mw^s,\xi}) =\operatorname{Vol} (S_{M,\xi})$. To check this, we use the polar
coordinates $x\rho+y=re^{i\theta}$ from \eqref{9.1}, leading to formula \eqref{9.2} for $\operatorname{Vol}
(S_{M,\xi})$. Note that the inequalities defining the volume in polar coordinates only depend on $\ell(M)$, $r$
and $\theta_M-2\theta$, with the restriction $\theta\in [0, \pi]$.
Since $\ell(M)$, $\theta_M$ only depend on $M\rho$, it follows that $\operatorname{Vol}
(S_{Mw^s,\xi})=\operatorname{Vol} (S_{M,\xi})$.

To show $\operatorname{Vol} (S_{w M,\xi})=\operatorname{Vol} (S_{M,\xi})$, let $\gamma_{x,y}\in \SL_2(\R) $ be any
matrix with lower row $(x,y)\ne (0,0)$; note that $j(\gamma_{x,y}, \rho):=x\rho+y=r e^{i\theta}$ in polar
coordinates. The transformation $(x,y,z)\mapsto (x',y',z)$ with $(x',y')$
defined by $\gamma_{x',y'}=\gamma_{x,y} w^{-1}$ has $x'\rho+y'=j(\gamma_{x,y}, \rho) j(w^{-1},
\rho)=re^{i(\theta-\pi/3)}$, so in polar coordinates it corresponds to $(r,\theta)\mapsto
(r'=r,\theta'=\theta-\frac{\pi}{3})$. Since $M \mapsto wM$ results in $\theta_M\mapsto \theta_M-\frac{2\pi}{3}$, and the
inequalities \eqref{9.3} defining $S_{M,\xi}$ involve only $\theta_M-2\theta$, it follows that
the transformation $(x,y,z)\mapsto (x',y',z)$ above
maps the volume $S_{M,\xi}$ onto a volume $S_{w M,\xi}^{\prime}$, defined like $S_{w M,\xi}$ but with the range
 $\theta\in [0,\pi]$ replaced by $\theta'\in [-\frac{\pi}{3},\frac{2\pi}{3}]$.  Since in the formula for $B_M(\xi,t)$ following
\eqref{9.3} the integrand has period $\pi$, we conclude $\operatorname{Vol} (S_{w M,\xi})=\operatorname{Vol}
(S_{M,\xi})$,

This concludes the proof of \eqref{1.9}. The formula in Theorem \ref{Thm1.1} follows from the results of Section
\ref{sec9}.

\appendix
\section{Arithmetic description of closed geodesics through $\rho$}\label{appendix}

In this appendix we discuss the connection between the hyperbolic lattice centered at $\rho$ and closed geodesics
on the modular surface passing through $\Pi(\rho)$ where $\Pi:\H\rightarrow \H/\Gamma$ is the projection map. In
the case of the hyperbolic lattice centered at $i$, the corresponding geodesics are the reciprocal geodesics
studied by Fricke and Klein \cite[Section 2]{BPPZ}. We  similarly describe the primitive closed geodesics passing
through $\Pi(\rho)$, which have an interesting arithmetic structure.

Closed geodesics on the modular surface correspond to conjugacy
classes $\{ g \}$ of hyperbolic elements $g\in \Gamma$. If $\a_g\subseteq \H$ is the axis of
$g$ (the semicircle connecting the two fixed points of $g$ on the real axis), then the
geodesic corresponding to $g$ on $X$ is $\Pi(z_0\rightarrow g z_0)$ for any fixed point
$z_0\in \a_g$. We are interested in geodesics passing through $\Pi(\rho)$. Let
$\RR$ denote the set of conjugacy classes of hyperbolic elements which contain a matrix $g$
whose axis passes through $\rho$. Let $\RRp\subseteq \RR$ be the subset of primitive
conjugacy classes. We will give an arithmetic description of $\RRp$.

Let $g=\left( \begin{smallmatrix} A & B \\ C & D \end{smallmatrix} \right)$ be a
primitive hyperbolic matrix whose axis $\a_g$ passes through $\rho$. The fixed points
$\lambda>\overline{\lambda}$ of $g$ satisfy the equation
\[
C\lambda^2+(D-A) \lambda -B=0 .
\]
Imposing the condition $(\lambda-\overline{\lambda})^2=|\rho-\lambda|^2+|
\overline{\rho}-\lambda|^2$ we conclude that
\begin{equation}\label{A1}
 \rho\in \a_g \iff D-A=2(B-C).
\end{equation}
The matrices $wgw^2$, and $w^2 g w$ are also primitive. Their axes are the same as
$\a_g$ rotated by $\pm \frac{2 \pi}{3}$ around $\rho$, hence the class
$\{g\}$ contains a matrix, still denoted by $g$, with $g\rho \in {\bf
I}$ (the first region in Section \ref{sec10}). This is equivalent to $\lambda \in
( \frac{1}{2} ,2)$ and $\operatorname{Re} (g\rho)> \frac{1}{2}$, which is further equivalent with $A,B,C,D$ being all
positive or all negative.

In conclusion each class $\{h\} \in \RRp$ contains a matrix $g$ as above with positive
entries satisfying \eqref{A1}, and so we are left to describe the set
of such matrices and determine when two such matrices are conjugate.

The condition \eqref{A1}, together with $AD-BC=1$, implies that $(A+D)^2-4(B^2+C^2-BC)=4$.
Writing $k=(B,C)$, $B=kB_0$, $C=k C_0$, $T:=A+D$, the pair $(T,k)$ is a
solution to Pell's equation
\begin{equation}\label{A2}
T^2-4k^2\Delta=4,
\end{equation}
with $\Delta=B_0^2+C_0^2-B_0 C_0$. In fact, $(T,k)$ is the minimal positive
solutions since $g$ is primitive. Direct computation using \eqref{A1} shows that
\begin{equation}\label{A3}
 \cosh d(\rho, g\rho)= \frac{T^2}{2}-1.
\end{equation}
We are led to define the set
\begin{equation}\label{A4}
\begin{split}
\DD_\rho &:=\bigg\{(B_0,C_0) : \begin{matrix} B_0,C_0>0, \ (B_0,C_0)=1,\\ \Delta=B_0^2+C_0^2-B_0 C_0
\text{ not a square } \end{matrix}\bigg\}
=\bigcup_{\Delta\in D_\rho}\DD_\Delta ,
\end{split}
\end{equation}
where $\DD_\Delta$ is the finite subset of pairs $(B_0, C_0)$ as above having fixed
$\Delta$. We denoted by $\DD_\rho$ the set of possible such
$\Delta$, which is the same as the set of positive numbers all of whose prime factors are
congruent to 1 mod 3, or the prime 3 appearing to the first power. The cardinality of $\DD_\Delta$ is
$2^{1+\nu}$ with $\nu$ the number of distinct prime factors $p\equiv 1 \pmod 3$ of
$\Delta$.
We conclude that there is a parametrization
\[
\varphi: \DD_\rho\rightarrow \RRp, \quad \varphi(B_0,C_0)=\left\{\left(
\begin{matrix} \frac{T}{2}-k(B_0-C_0) & k B_0 \\ k C_0  & \frac{T}{2}+k(B_0-C_0)
\end{matrix} \right)\right\} ,
\]
where $(T,k)$ is the smallest positive solution of Pell's equation $T^2-4k^2\Delta=4$ (the
minimality of $(T,k)$ ensures that the image of $\varphi$ consists of primitive conjugacy
classes only).

We are left to determine, for each primitive hyperbolic $g\in \Gamma$ satisfying
\eqref{3.1} and having positive entries, the set of $h\neq g$ with positive entries, having
$\{g\}=\{h\}$.
Assume therefore that $h=\gamma^{-1} g \gamma$. Then $\gamma$ maps $\a_h$ onto $\a_g$ (as
it can be seen by looking at what $\gamma$ does to the endpoints of $\a_h,\a_g$), and we
therefore have that $\gamma \rho\in \a_g$. Replacing $\gamma$ by $g^n \gamma$ for an
appropriate $n$, we can assume that $\gamma \rho \in (\rho\rightarrow
g\rho)$. Writing $g=\gamma \gamma'$, it follows that $h=\gamma'\gamma$. Since $\gamma$
maps $\a_h$ onto $\a_g$, and since the point $g\rho =\gamma\gamma'\rho\in(\gamma \rho
\rightarrow \gamma \gamma' \gamma \rho)$, it follows that $\gamma'\rho\in (\rho\rightarrow
h\rho)$.

Therefore the number of hyperbolic $h$ with $\{h\}=\{g\}$ and $h\rho \in {\bf I}$  is the
same as the number of points $\gamma \rho\in (\rho\rightarrow g\rho)$ (open geodesic
segment) with $\gamma\in\Gamma$ (compare with Lemma 2 of \cite{BPPZ}), that is the same as
the number of decompositions $g=\gamma \gamma'$ with $\gamma, \gamma'$ having positive
entries (the positivity follows from $\gamma \rho, \gamma'\rho \in {\bf I}$ and \eqref{10.1}). We will show that
there are 0, 1, or 3 such points, depending on arithmetic conditions on $\Delta$.

Let $X,Y,Z$ be the coordinates \eqref{2.1} for such a $\gamma\in \Gamma$, so $X,Y,2Z\in \N$.
We have
$$\gamma \rho \in \a_g \iff
\frac{2Z-X}{2Z-Y}=\frac{B}{C}=\frac{B_0}{C_0},
$$
therefore $2Z=X+u B_0$, $2Z=Y+u C_0$. The equation
$4XY-4Z^2=3$ becomes, after setting $t=6 Z-2u(B_0+C_0)$,
\begin{equation}\label{A5}
 t^2-4 u^2 \Delta=9.
\end{equation}
In terms of solutions $(t,u)$ with $t>0$ we find
\begin{equation}\label{A6}
\begin{split}
2Z & =\textstyle\frac{1}{3}\big(t+2u(B_0+C_0)\big),\quad X=\textstyle\frac{1}{3}\big(t+u (2C_0-B_0)\big), \\
Y & =\textstyle\frac{1}{3}\big(t+u(2B_0-C_0)\big).
\end{split}
\end{equation}
If $(3,u)=1$, the sign of $u$ is determined by the condition $t\equiv
2u(B_0+C_0)\pmod{3}$ which ensures that $X,Y,2Z \in\Z$. Notice that $X,Y>0$, so the triple
$(X,Y,Z)$ indeed determines a matrix $\gamma$ with $\gamma \rho\in \a_g$.  By
\eqref{2.2} we find
\begin{equation}\label{A7}
 \cosh d(\rho, \gamma \rho)=\frac{t}{3} .
\end{equation}
We distinguish two types of solutions, depending on whether 3 divides $u$ or not.

{\bf Case I:} 3 divides $u$. Letting $u=3u'$, $t=3t'$, with $(t',u')$ a solution of
$$t'^2 - 4 \Delta u'^2=1,$$ we have from \eqref{A6} that $2Z>X, 2Z>Y$ when $u'>0$, so the
point $\gamma\rho$ is on the same side of $\a_g$ as $g\rho$.
Since $d(\rho, \gamma \rho)=t'$, to determine when $\gamma \rho \in
(\rho\rightarrow g\rho)$ we distinguish two cases : if
the minimal positive solution $(T,k)$ of \eqref{A2} has $k$ even, then $2t'=T, 2u'=k$
and $d(\rho, g \rho)=2d(\rho, \gamma \rho)$, and we find that $\gamma \rho$ is the midpoint
of $(\rho\rightarrow g\rho)$. On the other hand if $k$ is odd, then the
minimal solution $(t',u')$ has $t'=\frac{T^2}{2}-1$, in which case $\gamma=g$ and there are no
points $\gamma \rho$ on $(\rho\rightarrow g\rho)$.

{\bf Case II:} $(3,u)=1$. In this case we necessarily have $(3,\Delta)=1$. Writing
\eqref{A5} as $N(\alpha)=9$ with $\alpha=t+2u\sqrt{\Delta}$, and assuming a
solution $\alpha_0=t_0+2u_0\sqrt{\Delta}$ exists with $(3,u_0)=1$, then all solutions are
$\alpha=\alpha_0 ( \frac{T}{2}+k\sqrt{\Delta})$, with $(T,k)$ solution of
\eqref{A2} with $k$ even. The corresponding points $\gamma \rho$ to these solutions lie on the axis $\a_g$ at
distance $\frac{1}{2} d(\rho, g\rho )$ apart if the minimal solution $(T,k)$  of
\eqref{A2} has $k$ even, or at distance $d(\rho, g\rho)$ apart if the minimal solution has $k$ odd. In the former
case we find two points $\gamma \rho$ on $(\rho\rightarrow g\rho)$, at distance $\frac{1}{2} d(\rho, g\rho)$ apart, while
in the latter only one. These are distinct from the point found in Case I, since here $\cosh d(\rho, \gamma \rho)$ is
not integral by \eqref{A7}.

The other solutions come from other generators of the ideal $\pp_3^2$, and geometrically
the corresponding points $\{\gamma \rho\}$ are the translates $\{g^n \gamma_0
\rho\}_{n\in\Z}$ along $\a_g$ for a fixed $\gamma_0$. Therefore there is a point
$\gamma \rho \in (\rho, g\rho)$, and this point is distinct from the
midpoint found in {\bf Case I}, since $\cosh d(\rho, \gamma \rho)\not\in \Z$.

We conclude that there is a partition of the set $D_\rho$ in \eqref{A4}:
\[
 D_\rho=D_\rho^{0,0}\cup D_\rho^{1,0}\cup D_\rho^{0,1}\cup D_\rho^{1,1} ,
\]
where $D_{\rho}^{0,\epsilon}$, respectively $D_{\rho}^{1,\epsilon}$ is the subset of
$\Delta$ for which the minimal positive solution $(T,k)$ of \eqref{A3} has $k$ even,
respectively odd, and $D_{\rho}^{\epsilon,0}$, respectively $D_{\rho}^{\epsilon,1}$, is
the subset of $\Delta$ such that \eqref{A5} has a solution with $(3,u)=1$, respectively it
does not have such a solution.  From the preceding discussion we conclude that the
restriction of the parametrization
\[
\varphi:\DD_\Delta \rightarrow \RRp_\Delta
\]
is 1-1 if $\Delta\in D_\rho^{1,1}$ (no lattice points on $(\rho\rightarrow g \rho)$), 2-1
if $\Delta\in D_\rho^{1,0}$  or $\Delta\in D_\rho^{0,1}$ (one point on $(\rho\rightarrow g
\rho)$),  and 4-1 if $\Delta\in D_\rho^{0,0}$ (three points on $(\rho\rightarrow g
\rho)$). We denoted by $\RRp_\Delta$ the image of $\varphi$ restricted to
$\DD_\Delta$.

\noindent{\bf Examples. I.} $\Delta=3$, $\DD_\Delta=\{(2,1), (1,2)\}$. The minimal
solution of $T^2-12k^2=4$ is $(T,k)=(4,1)$ with $k$ odd and $3|\Delta$ so $\Delta\in
D_{\rho}^{1,1}$.
Therefore there are two conjugacy classes in $\RRp_3$ with representatives:
$\varphi(2,1)=\left( \begin{smallmatrix} 1 &
2 \\ 1 & 3 \end{smallmatrix} \right)$, $\varphi(1,2)=\left( \begin{smallmatrix} 3 & 1
\\ 2 & 1 \end{smallmatrix} \right).$

\noindent{\bf II.} $\Delta=7$, $\DD_\Delta=\{(3,1), (3,2), (1,3), (2,3)\} $.   The minimal
solution of $T^2-28 k^2=4$ is $(T,k)=(16,3)$, and the equation $t^2-28 u^2 =9 $ has
solution $(t,u)=(11,\pm 2)$, so $\Delta\in D_\rho^{1,0}$. From Case II we find
$g=\varphi(3,1)=\gamma \gamma',\  h=\varphi(3,2)=\gamma' \gamma$, for
$\gamma=\left( \begin{smallmatrix} 1 & 2 \\ 1 & 3 \end{smallmatrix} \right)$, $\gamma'=\left(
\begin{smallmatrix} 1 & 4 \\ 1 & 5 \end{smallmatrix} \right)$,
so there are two conjugacy classes in $\RRp_\Delta$, $\{g\}=\{h\}$ and
$\{\widetilde{g}\}=\{\widetilde{h}\}$.

\noindent{\bf III.} $\Delta=21$, $\DD_\Delta=\{(5,1), (5,4), (1,5), (4,5)\} $. The minimal
solution of $T^2-84 k^2=4$ is $(T,k)=(2\cdot 55,12)$ and $3|\Delta$, so $\Delta\in
D_\rho^{0,1}$. From Case I we find
$g=\varphi(5,1)=\gamma \gamma'$, $h=\varphi(3,2)=\gamma' \gamma$, for
$\gamma=\left( \begin{smallmatrix} 3 & 4 \\ 5 & 7 \end{smallmatrix} \right)$, $\gamma'=\left(
\begin{smallmatrix} 1 & 8 \\ 1 & 9 \end{smallmatrix} \right)$,
so there are two conjugacy classes in $\RRp_\Delta$, $\{g\}=\{h\}$ and
$\{\widetilde{g}\}=\{\widetilde{h}\}$.

\noindent{\bf IV.} $\Delta=13$, $\DD_\Delta=\{(4,1), (4,3), (1,4), (3,4)\}$.
The minimal solution of $T^2-52 k^2=4$ is $(T,k)=(2\cdot 649,180)$, and the equation $t^2-52 u^2
=9 $ has solution $(t,u)=(29,\pm 4)$, so $\Delta\in D_\rho^{0,0}$. We have
$\varphi(4,1)=\left(\begin{smallmatrix} 109 & 720 \\ 180 & 1189
\end{smallmatrix}\right)=\gamma_1 \widetilde{\gamma}_1 =\gamma_2 \gamma_2'=\gamma_3 \gamma_3'$,
with $\gamma_1=\left(\begin{smallmatrix} 20 & 3 \\ 33 & 5 \end{smallmatrix}\right)$,
$\gamma_2=\left(\begin{smallmatrix} 2 & 1 \\ 3 & 2 \end{smallmatrix}\right)$,
$\gamma_2'=\left(\begin{smallmatrix} 38 & 251 \\ 33 & 218 \end{smallmatrix}\right)$,
$\gamma_3=\left(\begin{smallmatrix} 43 & 66 \\ 71 & 109 \end{smallmatrix}\right)$,
$\gamma_3'=\left(\begin{smallmatrix} 1 & 6 \\ 1 & 7 \end{smallmatrix}\right)$. The first
decomposition comes from Case I, while the last two from Case II.
We have $\varphi(1,4)= \widetilde{\gamma}_1 \gamma_1, \varphi(3,4)=\gamma_2'
\gamma_2, \varphi(4,3)=\gamma_3' \gamma_3$, so $\RRp_{13}$ contains one
conjugacy class.

In Case I, the midpoint of $(\rho\rightarrow g\rho)$ can be determined as follows.
Let $g=\left(\begin{smallmatrix} A & B \\ C & D \end{smallmatrix}\right)$
be a hyperbolic matrix whose axis passes through $\rho$, having positive entries. Then a
matrix $\gamma\in \SL_2(\R)$ such that $\gamma \rho$ is the midpoint of $(\rho\rightarrow g
\rho)$ has coordinates $(X,Y,Z)$ given by \eqref{2.1} as
$X=A+\frac{B}{2}$, $Y=D+\frac{C}{2}$, $Z=\frac{A}{2}+B=\frac{D}{2} +C .$

\section*{Acknowledgments}
The first author was partly supported by CNCS - UEFISCDI grant PN-II-RU-TE-2011-3-0259.
The second author was
partly supported by European Community grant PIRG05-GA-2009-248569 and CNCS - UEFISCDI grant
PN-II-RU-TE-2011-3-0259.

\end{document}